\documentclass[11pt,english]{article}

\usepackage{hyperref}
\hypersetup{
    unicode      = {false},
    pdftoolbar   = {true},
    pdfmenubar   = {true},
    pdffitwindow = {true},
    pdftitle     = {IDS},
    pdfauthor    = {Xavier Lebeuf},
    pdfsubject   = {IDS},
    pdfkeywords  = {IDS},
    pdfnewwindow = {true},
    colorlinks   = {true},
    linkcolor    = {blue},
    citecolor    = {blue},
    filecolor    = {black},
    urlcolor     = {blue},
    breaklinks   = {true}
}

\usepackage{amsmath}
\usepackage{amsthm}
\usepackage{amssymb}
\usepackage{amstext}
\usepackage{multirow}
\usepackage{epsfig}
\usepackage{soul}
\usepackage{colortbl}
\usepackage{curves}
\usepackage{epic}
\usepackage{pgf}
\usepackage{geometry}
\usepackage{times}
\usepackage{float}
\usepackage{microtype}
\usepackage{cases}

\definecolor{Red}{rgb}{1,0,0}
\definecolor{Green}{rgb}{0,.6,0}
\definecolor{Blue}{rgb}{0,0,1}
\definecolor{Fushia}{rgb}{0.7,0,0.7}

\usepackage{xspace}                                             
\usepackage{array}                                              
\usepackage{bbm}                                                
\usepackage{subcaption}                                         
\usepackage{tikz}                                               
\usetikzlibrary{arrows, shapes.geometric, positioning, calc}    
\usetikzlibrary{decorations.markings}                           
\usepackage[short]{optidef}                                     
\usepackage{chngcntr}                                           
\usepackage[capitalise,noabbrev,nameinlink]{cleveref}

\tikzstyle{block}       = [rectangle, draw, minimum width=1.2cm, minimum height=0.7cm, text centered, align=center, fill=white]
\tikzstyle{extblock}    = [rectangle, draw, anchor=south west, minimum height=1cm, text centered, align=center, fill=gray!10]
\tikzstyle{texte}       = [text centered, align=center]
\tikzstyle{decision}    = [diamond, draw, aspect=1, text centered, align=center, fill=white]
\tikzstyle{arrow}       = [->, thick]
\tikzstyle{darrow}      = [<->, thick]

\newcolumntype{C}{>{$}c<{$}}                                    
\newcolumntype{L}{>{$}l<{$}}                                    
\newcolumntype{R}{>{$}r<{$}}                                    
\def\ic{{\hat{\imath}}}                                         
\def\jc{{\hat{\jmath}}}                                         
\def\1{{\mathbbm{1}}}                                           
\def\R{{\mathbb{R}}}                                            
\def\N{{\mathbb{N}}}                                            
\def\RU{\overline{\R}}                                          
\newcommand{\solar}[1]{${\sf solar #1}$\xspace}                 
\newcommand{\nomad}{{\sf NOMAD}\xspace}                         
\newcommand{\inter}[1]{{\sf Inter\nobreakdash-#1}\xspace}                   
\newcommand{\mads}{{MADS}\xspace}                               
\newcommand{\ids}{{\sf IDS}\xspace}                             
\newcommand{\dids}{{\sf DIDS}\xspace}                           
\newcommand{\refP}[1]{\hyperlink{prob:p}{$\mathcal{P}(#1)$}\xspace}                  
\newcommand{\refPbar}{\hyperlink{prob:approxP}{$\overline{\mathcal{P}}$}\xspace}     
\newcommand{\refQ}{\hyperlink{prob:Q}{$\mathcal{Q}$}\xspace}                         

\newcounter{algocounter}                                                
\counterwithin{algocounter}{section}                                    
\renewcommand{\thealgocounter}{\thesection.\arabic{algocounter}}        
\newenvironment{algo}{\refstepcounter{algocounter}}{}                   
\newtheorem{assumption}{Assumption}                             
\newtheorem{theorem}{Theorem}
\newtheorem{lemma}[theorem]{Lemma}

\newtheorem{definition}{Definition}

\usepackage[textwidth=27mm, backgroundcolor=yellow, linecolor=orange, textsize=tiny]{todonotes}

\title{{
Multi-fidelity constraints in blackbox optimization \\
}
}

\author{
    \href{mailto:alarie.stephane@hydroquebec.com}{St\'ephane Alarie}\thanks{
        \href{https://www.hydroquebec.com/about/}{Hydro-Qu\'ebec}, \href{mailto:alarie.stephane@hydroquebec.com}{\tt alarie.stephane@hydroquebec.com}}
    \and
    \href{mailto:Charles.Audet@gerad.ca}{Charles Audet}\thanks{
        \href{https://www.gerad.ca}{GERAD} and \href{https://www.polymtl.ca}{Polytechnique Montr\'eal},
        \href{https://www.gerad.ca/Charles.Audet}{\tt www.gerad.ca/Charles.Audet}}
    \and
    \href{mailto:diagomartinez.miguel@hydroquebec.com}{Miguel Diago}\thanks{
        \href{https://www.hydroquebec.com/about/}{Hydro-Qu\'ebec}, \href{mailto:diagomartinez.miguel3@hydroquebec.com}{\tt diagomartinez.miguel3@hydroquebec.com}}
    \and
    \href{mailto:Sebastien.Le.Digabel@gerad.ca}{S\'ebastien {Le~Digabel}}\thanks{
        \href{https://www.gerad.ca}{GERAD} and \href{https://www.polymtl.ca}{Polytechnique Montr\'eal},
        \href{https://www.gerad.ca/Sebastien.Le.Digabel}
        {\tt https://www.gerad.ca/Sebastien.Le.Digabel}}
    \and
    \href{mailto:xavier.lebeuf@polymtl.ca}{Xavier Lebeuf}\thanks{
        \href{https://www.gerad.ca}{GERAD} and \href{https://www.polymtl.ca}{Polytechnique Montr\'eal}, \href{mailto:xavier.lebeuf@polymtl.ca}{\tt xavier.lebeuf@polymtl.ca}
        {(corresponding author)}}
}

\begin{document}
\maketitle
\thispagestyle{empty}

\thispagestyle{empty}

\noindent
{\bf Abstract:}
{
This work studies constrained blackbox optimization problems that cannot be solved in reasonable time due to prohibitive computational costs.
This challenge is especially prevalent in industrial applications, where blackbox evaluations are costly.
However, constraints can be evaluated at various fidelities at a lower computational cost.
More specifically, this work targets situations in which the infeasibility of each individual constraint can be detected at lower fidelities, and where a large discrete number of fidelities are available.
Moreover, highly discontinuous problems which may fail to evaluate are considered, such that direct search methods are preferred to model-based ones.
To this effect, the Interruptible Direct Search (\ids) and the Dynamic Interruptible Direct Search (\dids) algorithms are proposed to leverage feasibility assessments from various fidelity levels to avoid high cost evaluations.
The results show highly increased performances from \nomad when it is paired with \ids or \dids.
}



\setcounter{page}{1}

\section{Introduction}
\label{sec:intro}

This work studies constrained multi-fidelity optimization problems of the form
\begin{flushright}
\begin{minipage}{0.15\linewidth}
    \hypertarget{prob:P}{$\mathcal{P}(\phi)$}
\end{minipage}
\begin{minipage}{0.8\linewidth}
    \begin{mini*}
        {\text{\footnotesize $x\in \Omega$}}{f(x,\phi),\qquad\text{ where }\quad\Omega=\{x\in X:c_j(x,\phi)\leq 0, j\in J\},}
        {}{}
    \end{mini*}
\end{minipage}
\end{flushright}
in which $X\subseteq \R^n$ is the bound constrained domain of the objective function~$f:X\times[0,1]\rightarrow \RU=\R\cup\{\infty\}$ and of the $m$ relaxable quantifiable~\cite{LedWild2015} constraint functions~$c_j:X\times[0,1]\rightarrow \RU$, for~$j\in J=\{1,2,\dots,m\}$.
The set of feasible points is denoted by~$\Omega$.
The parameter~$\phi \in [0,1]$ is the \emph{fidelity level} at which these functions are evaluated.
Hence, by convention, the problem that one wishes to solve is denoted by~\refP{1} and is called the \emph{truth}.
The objective and constraint functions are provided by a blackbox process.
They have no accessible analytical formulation, they might be highly discontinuous, and their derivatives are unavailable or may be non-existent~\cite{AuHa2017}.
This is why a direct search approach is adopted~\cite{DzRiRoZe2025}.
Moreover, the blackbox is expensive to evaluate, and may fail to execute.
The set~$\RU$ allows the use of common blackbox optimization tools.
Notably,~$f$ is assigned an infinite value by an extreme barrier method at infeasible points, and the vector~$c(x,\phi)=(c_1(x,\phi),c_2(x,\phi),\dots,c_m(x,\phi))$ may posses one or more infinite values when an evaluation fails.

The specificity of this work is that Problem~\refP{\phi} is {\em multi-fidelity}, meaning that evaluating the blackbox requires not only specifying a trial point~$x \in X$,
but also selecting a fidelity value~$\phi \in [0,1]$ that controls the accuracy and computational cost of the evaluation.
Lower fidelities correspond to lower precision, and generally, lower evaluation cost, and vice-versa.
An evaluation using fidelity~$\phi<1$ can be interpreted as calling a static surrogate model~\cite[Ch.~13]{AuHa2017} that provides point-wise approximations of~$f$ and~$c$.
The computational cost required to evaluate a trial point~$x\in X$ using fidelity level~$\phi\in[0,1]$ is denoted by~$\lambda(x,\phi)$, where~$\lambda:X\times[0,1] \rightarrow \R_+$.
This cost function is assumed to be unknown, but is usually increasing with respect to $\phi$.

This work considers a finite and discrete subset of~$L\in\N$ fidelities.
It is described by the set of fidelity indices~$I = \{1,2,\dots,L\}$, where the sequence~$\{\phi_i\}_{i\in I}$ is strictly increasing and~$\phi_L=1$.
Two novel optimization algorithms are introduced: Interruptible Direct Search (\ids) and Dynamic Interruptible Direct Search (\dids).
They propose a strategy to reduce the computational cost of solving the optimization problem by identifying and exploiting the minimal fidelities required to determine that a trial point is deemed infeasible.
This approach avoids costly high-fidelity evaluations in a context where solving the true Problem~\refP{1} directly is impossible, as the computational effort is prohibitive.
While \ids is applicable to any multi-fidelity problem, \dids targets problems with {\em intermediary outputs}.

\begin{definition}
    A multi-fidelity blackbox with domain~$X$ is said to return {\em intermediary outputs} for~$I$ if during an evaluation of trial point~$x\in X$ at fidelity~$1$, all~$f(x,\phi_i)$ and~$c(x,\phi_i)$ values for~$i\in I$ become sequentially available during an evaluation. The evaluation can be interrupted after reaching~$\phi_1$,~$\phi_2$ up until the maximal fidelity~$\phi_L=1$.
    \label{def:intermediaire}
\end{definition}
Such a blackbox can easily be built for a stochastic problem where the fidelity controls the number of Monte-Carlo (MC) draws using Sample Average Approximation (SAA)~\cite{HoBa2014}.
When one queries a stochastic blackbox function~$f(x)$ at a point~$x\in X$ for a MC draw, a noisy value~$f_{\xi}(x)$ is obtained, where~$\xi$ is a random variable.
A set of~$L$ increasing amounts of MC draws~$\{\eta_i\}_{i\in I}$ that correspond to~$L$ fidelities is determined, where~$\eta_L$ is considered as the truth and~$\phi_i=\frac{\eta_i}{\eta_L}\in[0,1]$ for each~$i\in I$.
When using SAA, a predefined set of noise observations~$\{\xi_\omega\}_{\omega\in\{1,2,\dots,\eta_L\}}$ randomly sampled from the distribution of~$\xi$ is considered. This allows the definition of a deterministic multi-fidelity value for each~$x\in X$ and each~$i\in I$.
\begin{alignat}{2}
    f(x,\phi_i)=
    \begin{cases}
        \frac{1}{\eta_1} \sum\limits_{\omega=1}^{\eta_1}f_{\xi_\omega}(x) & \text{if }i=1, \\
        \frac{1}{\eta_i} \left( \eta_{i-1}f(x,\phi_{i-1}) + \sum\limits_{\omega=\eta_i-\eta_{i-1}}^{\eta_i}f_{\xi_\omega}(x) \right) & \text{otherwise.}
    \end{cases}
    \label{eq:saa}
\end{alignat}
The same process applies to~$c_j(x,\phi)$ for each~$j\in J$. A non-deterministic framework is presented in~\cite{AlAuBoLed2019}, where the observations~$\omega$ are not predefined.
When this new multi-fidelity blackbox is called at a trial point~$x\in X$, the first~$\eta_1$ draws are performed and~$f(x,\phi_1)$ and~$c(x,\phi_1)$ are made available.
If the evaluation is not interrupted, only the~$\eta_{2}-\eta_1$ next draws are performed, and the outputs for~$\phi_2$ are computed using~$f(x,\phi_1)$,~$c(x,\phi_1)$ and these new draws.
Then, to continue the evaluation, the next~$\eta_3-\eta_2$ draws are performed, and so on until~$\phi_L$ is reached.
This behaviour corresponds to \cref{def:intermediaire}.

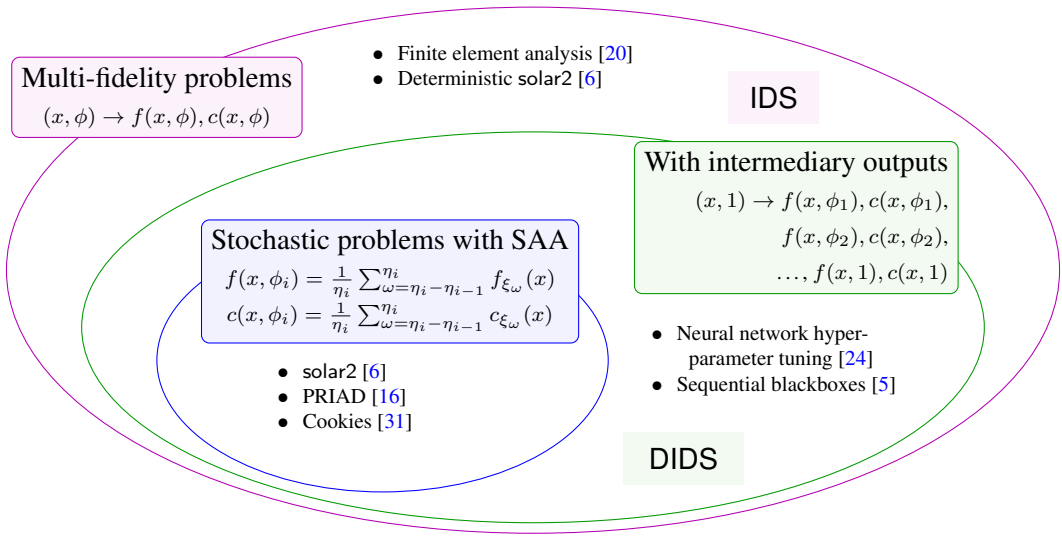
\begin{figure}[htb!]
    \begin{center}
    \begin{tikzpicture}[scale=1, transform shape]
    
        \draw[color=Blue]   (-2,-1.5) ellipse (3 and 1.75);
        \draw[color=Green]  (0,-1.06) ellipse (6 and 2.6);
        \draw[color=Fushia] (0,-0.3)  ellipse (7 and 3.5);
        
        \node[block,draw=Fushia, fill=Fushia!5!white,rounded corners=3pt]           (multifi)    at (-5,2)       {Multi-fidelity problems \\ 
            \scriptsize $(x,\phi) \rightarrow f(x,\phi),c(x,\phi)$};
        \node[block,align=right,draw=Green, fill=Green!5!white,rounded corners=3pt] (intermed)   at (3.5,0.4)    {With intermediary outputs \\ 
            \scriptsize $(x,1) \rightarrow f(x,\phi_1),c(x,\phi_1)$, \\ 
            \scriptsize $f(x,\phi_2),c(x,\phi_2)$, \\ 
            \scriptsize \dots, $f(x,1),c(x,1)$};
        \node[block,draw=Blue, fill=Blue!5!white,rounded corners=3pt]               (stochastic) at (-1.9,-0.45) {Stochastic problems with SAA \\     
            \scriptsize $f(x,\phi_i)=\frac{1}{\eta_i} \sum_{\omega=\eta_i-\eta_{i-1}}^{\eta_i}f_{\xi_\omega}(x)$ \\ 
            \scriptsize $c(x,\phi_i)=\frac{1}{\eta_i} \sum_{\omega=\eta_i-\eta_{i-1}}^{\eta_i}c_{\xi_\omega}(x)$};
        \node[block, minimum width=1.2cm, draw=none, fill=Fushia!5!white]            (ids)        at (3.2,2)      {\ids};
        \node[block, minimum width=1.6cm, draw=none, fill=Green!5!white]             (dids)       at (2,-2.8)     {\dids};

        \node[texte] at (-0.4, 2.4) {\scriptsize \begin{tabular}{l}
                $\bullet$ \, Finite element analysis~\cite{GuHaWaHuXi2020}\\
                $\bullet$ \, Deterministic \solar{2}~\cite{solar_paper}
            \end{tabular}};
        \node[texte] at (3.2, -1.5) {\scriptsize \begin{tabular}{l}
                $\bullet$ \, Neural network hyper- \\
                \hspace{0.5cm}parameter tuning~\cite{LaLed21} \\
                $\bullet$ \, Sequential blackboxes~\cite{G-2021-65}
            \end{tabular}};
        \node[texte] at (-2.5, -2) {\scriptsize \begin{tabular}{l}
                $\bullet$ \, \solar{2}~\cite{solar_paper}\\
                $\bullet$ \, PRIAD~\cite{DiLeMeCoRa2025}\\
                $\bullet$ \, Cookies~\cite{KGKSMS2017}
            \end{tabular}};
    \end{tikzpicture}
    \end{center}
    \caption{Classification of a few multi-fidelity problems with domains of application of \ids and \dids.}
    \label{fig:classification}
\end{figure}

\cref{fig:classification} illustrates a classification of some known problems with respect to \cref{def:intermediaire} and the applicability of SAA, as well as the classes for which the \ids and \dids algorithms are applicable.
PRIAD is an electrical equipment maintenance optimization problem, and \solar{} is a concentrated solar power plant design optimization problem.
\begin{itemize}
    \item \textbf{Multi-fidelity problems.}
    In finite element analysis, fidelity usually controls the coarseness of the mesh~\cite{GuHaWaHuXi2020}.
    The \solar{2} problem~\cite{solar_paper} is deterministic when its seed and number of MC draws are predetermined.
    Then, the fidelity controls the convergence criteria of numerical methods.
    \item \textbf{With intermediary outputs.}
    A neural network hyper-parameter optimization problem can be formulated such that the accuracy is returned after different number of epochs (different fidelities) during training~\cite{LaLed21}.
    Sequential blackboxes are problems where each output is given by a distinct blackbox which corresponds to a fidelity level, and they are each called sequentially~\cite{G-2021-65}.
    \item \textbf{Stochastic problems.}
    PRIAD's blackbox problem~\cite{DiLeMeCoRa2025}, \solar{2}~\cite{solar_paper} and the cookie recipe optimization problem~\cite{KGKSMS2017} are stochastic multi-fidelity problems where SAA may be applied.
    Additionally, this last problem is a laboratory experiment where intermediary outputs are easily accessible during a blackbox evaluation.
\end{itemize}

\subsection{Motivation}

The primary motivation for this work is an asset management blackbox optimization problem encountered at Hydro-Qu\'ebec as part of the PRIAD project~\cite{PRIAD_CoBLALDeKoMe2020,GaChKoCoHeBlDeAb2021,PRIAD_KoMeCoGaVoAlDeBl2021,DiLeMeCoRa2025}.
Its objective is to develop a periodic maintenance strategy optimization framework that relies on blackbox optimization methods.
As the project is under development, the blackbox simulator involved is not yet available.
Preliminary tests suggest that a single high fidelity evaluation of this simulator could take up to 145 days~\cite{PRIAD_KoMeCoGaVoAlDeBl2021}.
For an optimization of 2,000 evaluations, parallel computing can reduce the optimization time to one week if 45,000 CPUs are used~\cite{PRIAD_KoMeCoGaVoAlDeBl2021}.
To this effect, this work presents a new cost reduction method based on the idea that the feasibility of some constraints may be estimated with low fidelities.
When it is estimated that a point is infeasible from low fidelity information, high fidelity information is not computed to avoid high evaluation costs.

The second motivation of this work is to improve on the \inter{DS} algorithm presented in~\cite{AlAuDiLedLe23}.
\inter{DS} served as a first step towards an algorithmic approach that exploits information from a broad range of available fidelities in the context of direct search methods for constrained blackbox optimization.
The \ids algorithm of the present work is a strict upgrade from \inter{DS} in a theoretical analysis sense.
Many assumptions necessary to~\cite{AlAuDiLedLe23} are lifted for \ids, and in addition, some pathological cases are avoided without negatively impacting the algorithmic performance.
In \cref{fig:classification}, \inter{DS} would lie in the same class as \ids.
An important observation from~\cite{AlAuDiLedLe23} is that the sample point selection, required to solve an assignment subproblem, is crucial to the method.
A second algorithm, \dids, is proposed to periodically solve the subproblem as the optimization process proceeds.
This new method is expected to perform better, but it requires the blackbox problem to have intermediary outputs, as described in \cref{def:intermediaire}.

\subsection{Contribution}

Scientific literature on the subject of costly multi-fidelity blackbox problems predominantly studies the unconstrained case, or considers constraints along with a penalty in the objective function.
Moreover, multi-fidelity frameworks sometimes allow for more than two fidelity levels, but a single high fidelity and a single low fidelity are almost always considered in benchmarks.
Lastly, to the knowledge of the authors, all current multi-fidelity approaches use a low fidelity source to approximate gradients. This is usually achieved by fitting a model on the low fidelity data~\cite{lili2024}. 
Conversely, this research approaches multi-fidelity problems under radically different lens:
\begin{itemize}
    \item A finite and discrete set of fidelities is considered. The benchmarks in \cref{sec:results} use up to 11 fidelity levels.
    \item Direct-search methods are preferred to model-based approaches due to noisy, highly discontinuous and unpredictable blackbox problems which may fail to evaluate.
    \item Constraints are handled directly rather than being penalized in the objective function. As a matter of fact, the handling of constraints is the main focus of this work, and multi-fidelity information of the objective function is not used to improve the optimization process.
\end{itemize}

Future work will integrate multi-fidelity information of the objective function into the presented methods.

\subsection{Organization}

The document is structured as follows. 
\cref{sec:rev_lit} contains a short literature review on interruptible methods in blackbox optimization and on the use of multi-fidelity in the constrained case. \cref{sec:components} presents each of the new algorithmic components.
They are assembled in \cref{sec:inter_algos} to create the two distinct optimization algorithms, \ids and \dids.
Finally, \cref{sec:results} compares \ids and \dids to other methods on problems from the \solar{} benchmarking collection~\cite{solar_paper}.
On average, \dids outperforms \ids which in turn outperforms \inter{DS}.
Concluding remarks follow in the final section.

\section{Literature review}
\label{sec:rev_lit}

This work uses the KARQ (Known/hidden, A~priori/simulated, Relaxable/unrelaxable, Quantifiable/nonquantifiable) taxonomy of constraints~\cite{LedWild2015}.
Notably, a constraint is said to be a~priori if it has an explicit analytical formulation.
When evaluating a point, a~priori constraints are checked first at a negligible cost, and the blackbox process is only launched if they are satisfied.
When counting evaluations during an optimization, points where the blackbox process is not launched are not considered.
This work also considers direct search methods for blackbox optimization, which are categorized into directional,
mesh-based and line-search algorithms~\cite{DzRiRoZe2025}.
For these methods, the quantitative constraints~$c_j(x)\leq0,j\in J$ are often handled through the constraint violation function~$h:\R^m\rightarrow\RU$ introduced in~\cite{AuDe09a}:
\begin{equation*}
    h(x):=
    \begin{cases}
        \sum\limits_{j=1}^m(\max\{c_j(x),0\})^2 & \text{ if }x\in X \\
        \infty & \text{ otherwise.}
    \end{cases}
\end{equation*}

The extreme barrier (EB)~\cite{AuHa2017} algorithm is a two-phase method that first minimizes $h(x)$ in $X$ for as long as the incumbent solution is infeasible. 
When it reaches feasibility, the second phase minimizes the extended-value function
\begin{equation*}
    f_\Omega(x):=
    \begin{cases}
        f(x) & \text{ if }x\in \Omega \\
        \infty & \text{ otherwise}
    \end{cases}
\end{equation*}
on $X$ from the feasible starting point.
This approach is employed to reject infeasible points.
The progressive barrier~\cite{AuDe09a} is a more advanced constraint handling method.
It consists of maintaining two incumbent solutions, a feasible and an infeasible one.
The infeasible solution is progressively pushed towards the feasible domain by decreasing a threshold on~$h(x)$, above which trial point~$x\in X$ is rejected.
The infeasible incumbent often has a better objective function value, which leads to the discovery of new feasible incumbent solutions.

The two phase interruptible EB~\cite{G-2021-65} is a barrier method that applies to sequential blackboxes problems described above \cref{fig:classification}.
During an evaluation of trial point~$x\in X$, the constraint violation function~$h(x)$ is updated after each constraint value is computed.
When this value is greater than~$h(x^k)$ where~$x^k$ is the incumbent solution, the evaluation is interrupted.
In the case where the ordering of the blackboxes can be chosen, the hierarchical satisfiability EB~\cite{G-2021-65} suggests a sequence of optimization problems to solve with a similar interruption mechanism.
These methods are built on the principle that costly evaluations should be interrupted when it is recognized through intermediate processes that a solution will not benefit the search strategy~\cite{RaToMaThMaSe2010}.
This principle is relevant to multi-fidelity blackbox problems, since fidelity levels correspond to such intermediate processes.
The \inter{DS} algorithm~\cite{AlAuDiLedLe23} is designed to leverage multi-fidelity information to quickly identify infeasible points and interrupt their evaluation.
An important contribution introduced a trust region algorithm that incorporates multi-fidelity constraints and objective function into a sufficient decrease framework~\cite{March2012}.
An in-depth review of the relevant blackbox optimization and multi-fidelity literature is presented in~\cite{AlAuDiLedLe23}.

In the recent literature, a trust region method for stochastic problems where each level of fidelity corresponds to a different stochastic simulation subject to different numbers of MC draws is presented in~\cite{HaMu2025}.
Constrained stochastic problems are also tackled under a multi-fidelity framework by the Scout-Nd algorithm~\cite{AgRaKoBu2023}.
In the field of Bayesian optimization, the expected constrained improvement~\cite{GaKuXuWeCu2014} was introduced as an extension to the expected improvement for constrained problems.
It consists of multiplying the expected improvement with a probability of feasibility, and it was recently used in the context of multi-fidelity constrained Bayesian optimization~\cite{sMFBO2CoGP,KhBeBrDeFe2020,QuJiQiLeAn2024}.
Building on this work, new acquisition functions which are easier to implement by not solely relying on the expected constraint improvement were proposed~\cite{WaChHaPeSoPe2025}.

\section{Algorithmic components}
\label{sec:components}

The methods presented herein exploit multi-fidelity information from a blackbox problem with the primary goal of reducing evaluation costs during the optimization process. 
They are designed to work in conjunction with an existing direct search solver that does not natively support multi-fidelity evaluations. 
In order to efficiently leverage information across fidelity levels, it is essential that data from previous evaluations at various fidelities is accessible.  
\cref{sec:inter_algos} proposes the \ids and \dids algorithms, each applicable under different assumptions regarding the nature of this accessibility, and to different problems as classified in \cref{fig:classification}.
These algorithms share several components, which are described individually in the present section, before introducing the complete frameworks.

\cref{sec:fidelity_control} introduces the {\em fidelity controller algorithm}, which performs the feasibility estimations.
These estimations are guided by a key parameter of the fidelity controller: the {\em assignment vector}.
Then, \cref{sec:compute_a} describes the assignment vector computation process by solving an optimization subproblem.
\cref{sec:theorie_q} presents theoretical guarantees of the subproblem.

\subsection{The fidelity controller algorithm}
\label{sec:fidelity_control}

The fidelity controller algorithm is a wrapper around the blackbox.
Its purpose is to reduce evaluations costs by estimating the feasibility of a candidate using different fidelities with indices stored in~$I$, for the case when the solver does not natively handle fidelity.
Hence, rather than solving Problem~\refP{1} directly, the wrapper problem~\refPbar is provided to the solver.
Problem~\refPbar is an accelerated version of Problem~\refP{1}, with the intention of conserving feasibility.
Through this document, for any optimization problem~$\mathcal{P}$, \textit{evaluating~$\mathcal{P}$} is a shorthand used for \textit{evaluating the blackbox from which Problem~$\mathcal{P}$ is derived}.
\begin{flushright}
\begin{minipage}{0.15\linewidth}
    \hypertarget{prob:approxP}{$\overline{\mathcal{P}}$}
\end{minipage}
\begin{minipage}{0.8\linewidth}
    \begin{mini*}
        {x\in X}{\bar f(x)\quad\text{ s.t. }\quad \bar c_j(x)\leq 0, j\in J.}
        {}{}
    \end{mini*}
\end{minipage}
\end{flushright}
The functions~$\bar f(x)$ and~$\bar c(x)$ correspond to values of~$f(x,\phi)$ and~$c(x,\phi)$, where~$\phi\in[0,1]$ is determined by the fidelity controller.

Consider an assignment of the constraint~$c_j\leq0,j\in J$, to the fidelity level~$\phi_i,i\in I$, with vector~$a\in I^m$, where~$a_{j}=i$ if constraint~$c_j\leq0$ is assigned to~$\phi_i$.
An assignment of a constraint to a fidelity indicates the lowest fidelity for which it is estimated that the constraint's violation is accurately asserted.
\cref{sec:compute_a} shows how this assignment vector~$a$ is computed.
As~$I$ can be large, a subset of fidelities are selected, which depends on this vector~$a$.
This subset is denoted by~$\Phi(a)\subseteq\{\phi_i\}_{i\in I}$, the set of fidelities of interest for the fidelity controller. Two distinct descriptions of this subset and of~$I$ are given in \cref{sec:ids} for \ids and \cref{sec:dids} for \dids.

This assignment is at the core of the fidelity controller algorithm.
Whenever~\refPbar is evaluated at a point~$x\in X$ during the optimization,~\refP{\phi} is sequentially evaluated at increasing fidelity levels~$\phi\in\Phi(a)$ by the fidelity controller algorithm.
After an evaluation at fidelity~$\phi$, only the constraints assigned to~$\phi$ or a lower fidelity are checked.
If any of those constraints is violated, the sequence of evaluations is interrupted, and the evaluated point is deemed infeasible.
The goal is to identify infeasible points as cheaply as possible, and to stop investing computational costs into their evaluation.
When all evaluations of~\refP{\phi} are completed without interruption, the evaluated point is deemed feasible.

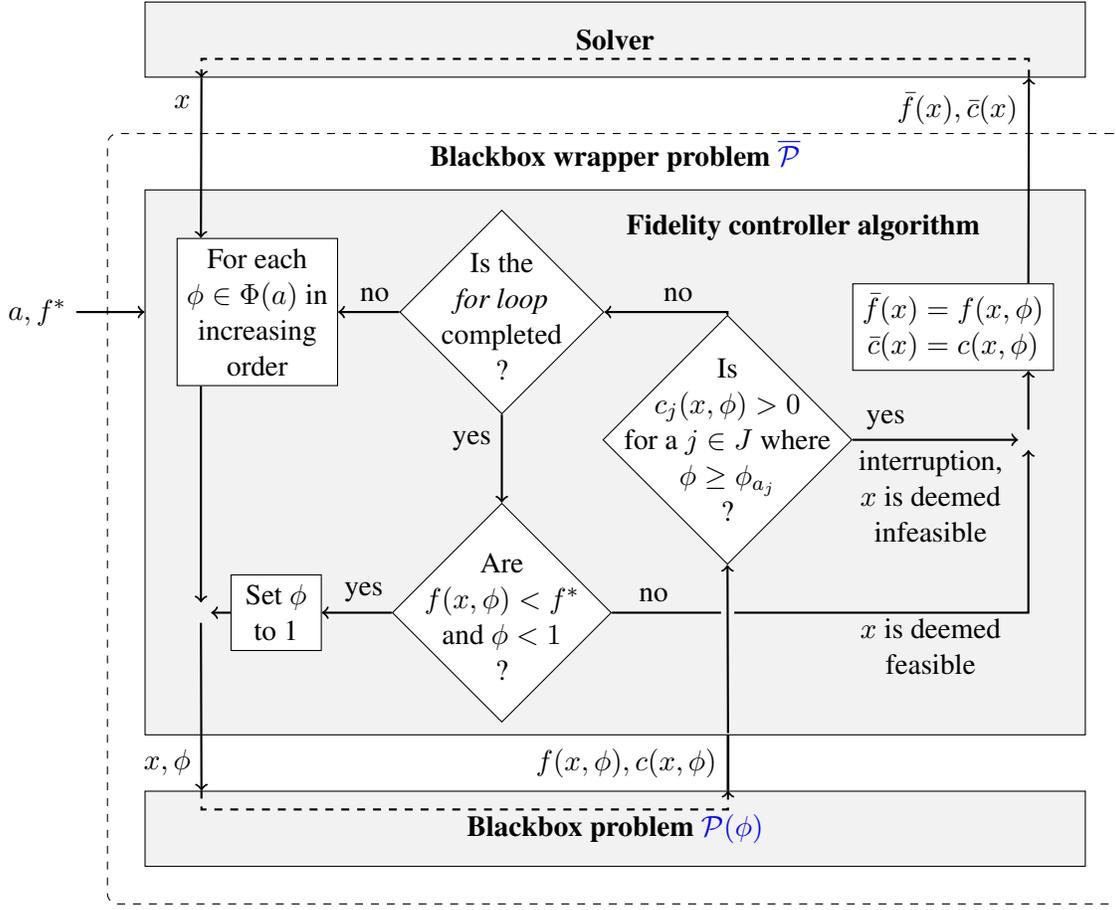
\begin{figure}[htb]
    \begin{center}
    \begin{tikzpicture}[scale=1, transform shape]
    
        \node[extblock, minimum width=12.5cm]                                                                                  (bb)     at (0,0)                   {\textbf{Blackbox problem }\refP{\phi}};
        \node[draw, dashed, inner sep=0.3cm, rounded corners, minimum width=13.5cm, minimum height=10.25cm, anchor=south west] (bbw)    at (-0.5,-0.5)             {};
        \node[extblock, minimum width=12.5cm, minimum height=7.25cm]                                                           (fico)   at (0,1.75)                {};
        \node[extblock, minimum width=12.5cm]                                                                                  (solver) at (0,10.5)                {\textbf{Solver}};
    
        \node[block]                      (forloop)     at ($(fico.west)+(1.5,2)$)     {For each \\ $\phi\in\Phi(a)$ in \\ increasing \\ order};
        \node[]                           (float)       at ($(forloop)-(0.75,4)$)      {};
        \node[decision, text width=2.4cm] (evalsupp)    at ($(float)+(4,0)$)           {};
        \node[decision, text width=2.8cm] (interrupter) at ($(forloop)+(6.25,-1.7)$)   {};
        \node[decision, text width=2.2cm] (loopend)     at ($(forloop)+(3.25,0)$)      {};
        \node[block]                      (phi1)        at ($(float)+(1,0)$)           {Set $\phi$ \\ to 1};
        \node[]                           (fin)         at ($(interrupter)+(4,0)$)     {};
        \node[texte]                      (inputs)      at ($(forloop.east)-(4,0)$)    {$a,f^*$};
        \node[block]                      (outputs)     at ($(fin)+(-1,1.5)$)          {$\bar f(x) = f(x,\phi)$ \\ $\bar c(x) = c(x,\phi)$};   
    
        \node[texte] at ($(bbw.north)-(0,0.3)$)      {\textbf{Blackbox wrapper problem }\refPbar};
        \node[texte] at ($(fico.north)+(2.5,-0.5)$ ) {\textbf{Fidelity controller algorithm}};
        \node[texte] at ($(loopend)-(0,0.05)$)       {Is the \\ \textit{for loop} \\ completed \\ ?};
        \node[texte] at ($(interrupter)-(0,0)$)      {Is \\ $c_j(x,\phi)>0$ \\ for a $j\in J$ where \\ $\phi\geq\phi_{a_j}$ \\ ?};
        \node[texte] at ($(evalsupp)-(0,0.05)$)      {Are \\ $f(x,\phi)<f^*$ \\ and $\phi<1$ \\ ?};
        \node[texte] at ($(fin)-(1.3,2.75)$)         {$x$ is deemed \\ feasible};
        \node[texte] at ($(fin)-(1.3,0.75)$)         {interruption, \\ $x$ is deemed \\ infeasible};
        
        \draw[arrow]  ($(forloop.south)-(0.75,0)$)  --                                            (float.north);
        \draw[arrow]  (loopend.south)               -- node[pos=0.3 ,left ]              {yes}    (evalsupp.north);
        \draw[arrow]  (loopend.west)                -- node[pos=0.4 ,above]              {no}     (forloop.east);
        \draw[arrow]  (interrupter.east)            -- node[pos=0.2 ,above,align=left ]  {yes}    (fin.west);
        \draw[arrow]  (evalsupp.east)               -| node[pos=0.05,above,align=left ]  {no}     (fin.south);
        \draw[arrow]  (evalsupp.west)               -- node[pos=0.4 ,above]              {yes}    (phi1.east);
        \draw[arrow]  (phi1.west)                   --                                            (float.east);
        \draw[arrow]  (interrupter.north)           |- node[pos=0.7 ,above]              {no}     (loopend.east);
        \draw[arrow]  (fin.north)                   --                                            ($(outputs.south)+(1,0)$);
        
        \draw[arrow]          ($(solver.south)+(-5.5,0)$)  -- node[pos=0.15,left ] {$x$}                     ($(forloop.north)-(0.75,0)$);
        \draw[arrow]          (float.south)                -- node[pos=0.83,left ] {$x,\phi$}                ($(bb.north)+(-5.5,0)$);
        \draw[arrow]          ($(outputs.north)+(1,0)$)    -- node[pos=0.85,left ] {$\bar f(x),\bar c(x)$}   ($(solver.south)+(5.5,0)$);
        \draw[arrow]          (inputs.east)                --                                                ($(fico.west)+(0,2)$);
        \draw[line width=1pt] ($(bb.north)+(1.5,0)$)       -- node[pos=0.5 ,left ] {$ f(x,\phi),c(x,\phi)$}  ($(fico.south)+(1.5,0)$);
        \draw[arrow, preaction={draw=gray!10, line width=5pt}] ($(fico.south)+(1.5,0.03)$) --                (interrupter.south);
        
        \draw[arrow, dashed] ($(bb.north)+(-5.5,0)$)    -- ($(bb.north)+(-5.5,-0.25)$)   -- ($(bb.north)+(1.5,-0.25)$)     -- ($(bb.north)+(1.5,0)$);
        \draw[arrow, dashed] ($(solver.south)+(5.5,0)$) -- ($(solver.south)+(5.5,0.25)$) -- ($(solver.south)+(-5.5,0.25)$) -- ($(solver.south)+(-5.5,0)$);
    
    \end{tikzpicture}
    \end{center}
    \caption{Flow chart diagram of the fidelity controller algorithm. The \textit{for loop} is completed when~$\phi$ is greater or equal to the greatest element of~$\Phi(a)$, where~$a\in I^m$ is the assignment vector.}
    \label{fig:fid_flow}
\end{figure}

When an interruption occurs, the most recent outputs are returned to the solver.
If the subset of fidelities~$\Phi(a)$ does not contain~$1$, then it is possible that a point deemed feasible becomes the new incumbent solution, while being infeasible in reality.
The assignment computation method proposed in \cref{sec:compute_a} is such that this event is rare, but it may occur.
Since finishing the optimization with an infeasible point as the incumbent is highly undesirable, an additional evaluation at~$\phi_L=1$ is performed whenever this situation could arise.
The complete process is described in Algorithm~\ref{algo:fico} and schematized in \cref{fig:fid_flow}.

\begin{algo}
\label{algo:fico}
\begin{center}
    \begin{table}[ht]
    \centering \renewcommand{\arraystretch}{1.2}
    \begin{tabular}{L}
        \hline
        \textbf{Algorithm \thealgocounter: }\text{Fidelity controller algorithm
        (Blackbox wrapper problem~\refPbar)} \\
        \hline
        \textbf{Inputs: } \text{trial point } x\in X; \text{ assignment vector } a\in I^m; \text{ incumbent value } f^*\in\RU \\
        \textbf{Outputs: } \bar f(x)\in \RU;\, \bar c(x)\in \RU^m \\
        \textbf{For each }\text{fidelity level }\phi\in \Phi(a)\text{ in increasing order} \\
        \left.\middle|
                    \begin{array}{l}
                        \textbf{Evaluate } f(x,\phi) \text{ and } c(x,\phi) \\
                        \textbf{If }\text{there exists a } j\in J \text{ such that }\phi\geq \phi_{a_j} \text{ and }c_j(x,\phi) > 0 \\
                        \left.\middle|
                                   \begin{array}{l}
                                        \textbf{Return }\bar f(x)=f(x,\phi),\,\bar c(x)=c(x,\phi)
                                   \end{array}
                        \right. \\
                    \end{array}
        \right. \\
        \textbf{If } \phi<1 \text{ and } f(x,\phi) < f^* \\
        \left.\middle|
            \begin{array}{l}
                \textbf{Evaluate } f(x,1) \text{ and } c(x,1)\\
            \end{array}
        \right. \\
    \textbf{Return }\bar f(x)=f(x,1),\,\bar c(x)=c(x,1) \\
    \hline
    \end{tabular}
    \end{table}
\end{center}
\end{algo}

This interruption mechanism is such that the solver may unknowingly receive blackbox output values from various fidelities.
Nevertheless, outputs for points deemed feasible are only returned to the solver at the highest fidelity level in~$\Phi(a)$.
It may help solver performance to use the extreme barrier function~$f_\Omega$, to reject points deemed infeasible, which are more susceptible to have outputs returned at low fidelity.
Note that the use of the extreme barrier with the fidelity controller algorithm allows for the use of some direct search unconstrained optimization solvers. 

The fidelity controller algorithm is presented with the intention of reducing computational costs by causing interruptions on infeasible points while avoiding costly evaluations at~$\phi_L=1$.
Constraints should be assigned to minimal fidelities, while being assigned to fidelities sufficiently high to perform feasibility estimations.
Moreover, it may occur that estimations are incorrect, which can have an impact on the solver's performance.
For these reasons, choosing a good assignment is not a trivial task.

\subsection{The assignment computation}
\label{sec:compute_a}

Before using the fidelity controller algorithm from \cref{sec:fidelity_control}, the assignment vector~$a\in I^m$ must be computed. This section proposes a method to compute this assignment that minimizes the computational cost of an evaluation of~\refPbar, subject to constraints on the probability that a feasibility estimation is erroneous.
Consider a finite set of sample points~$H$. The indicator function
\begin{equation*}
    \1(c_j(x,\phi)>0) :=
    \begin{cases}
        1\qquad\text{ if }c_j(x,\phi)>0 \\
        0\qquad\text{ otherwise}
    \end{cases} \qquad\forall j\in J
\end{equation*}

\noindent is used to define {\em representativity}, a useful concept to determine the lowest fidelity at which the feasibility of a constraint can be correctly identified.
\begin{definition}
    Fidelity~$\phi\in[0,1]$ is said to be representative for a constraint~$c_j(x,\phi)\leq0$,~$j\in J$ at point~$x\in X$ if
    \begin{equation*}
        \1(c_j(x,\phi_i)>0)=\1(c_j(x,1)>0)\quad\forall\,i\in I\text{ where }\phi_i\geq\phi.
    \end{equation*}
    \label{def:representativity}
\end{definition}
For example, if the constraint values of~$c_1$ for~$L=4$ fidelities are~$(c_1(x,\phi_i))_{i\in I}=(-0.2,10,0,-3)$, then~$\phi_3$ and~$\phi_4$ are the only two representative fidelities for constraint~$c_1\leq0$ at point~$x\in X$.

To compute the assignment vector~$a$, statistical estimations are performed using the sample set~$H\subset X$.
Define~$X^\text{ap}$ as the subset of points from~$X$ where all a~priori constraints are satisfied, and~$H^\text{ap}=H\cap X^\text{ap}$.
When~$H\cap\Omega=\varnothing$, all constraints are assigned to~$\phi_L=1$, i.e.,~$a_j=L$ for each~$j\in J$.
In that case,~\refPbar is identical to~\refP{1}.
If no information about feasible points is available, then multi-fidelity information as a whole should not be exploited.
This allows \dids and \ids to avoid cases that are pathological for \inter{DS}.
Otherwise, given Algorithm~\ref{algo:fico}'s interruption mechanism, the vector~$a\in I^m$ that minimizes the expected computational cost of an evaluation of~\refPbar is found by solving the assignment Subproblem~\refQ.

\begin{flushright}
\begin{minipage}{0.05\linewidth}
    \hypertarget{prob:Q}{$\mathcal{Q}$}
\end{minipage}
\begin{minipage}{0.9\linewidth}
    \begin{mini!}
        {a\in I^{m}}{f_\mathcal{Q}(a) = \sum_{i:\phi_i\in\Phi(a)}\left(\lambda_i\prod_{j:a_j<i}p_{a_jj}\right) \label{eq:objQ}}
        {}{}
        \addConstraint{a_{j}}{\geq i(j) \qquad}{\forall\,j\in J \label{eq:rep_is_1}}
    \end{mini!}
\end{minipage}
\end{flushright}
where, each~$p_{a_jj}$,~$j\in J$ is assumed independent, and for each~$i \in I$ and~$j \in J$,
\begin{alignat}{1}
    \lambda_i : & = \frac{1}{\lvert H^\text{ap}\rvert}\sum_{x\in H^\text{ap}}\lambda(x,\phi_i) \approx \mathbb{E}_{x\in X^\text{ap}}[\lambda(x,\phi_i)], \label{eq:kappa} \\
    p_{ij} : & = \frac{1}{\lvert H^\text{ap}\rvert}\lvert\{x\in H^\text{ap}:c_j(x,\phi_i)\leq0\}\rvert \approx \mathbb{P}[c_j(x,\phi_i)\leq0 \text{, } x\in X^\text{ap}], \label{eq:p} \\
    r_{ij} : & = \frac{1}{\lvert H\cap\Omega\rvert}\lvert\{x\in H\cap\Omega:\phi_i\text{ is representative for }c_j\leq0\}\rvert \label{eq:r} \\
    & \approx \mathbb{P}[\text{fidelity $\phi_i$ is representative for constraint $c_j\leq0$, }x\in \Omega], \nonumber \\
    i(j) & :=\text{min}\{i\in I:r_{ij}=1\}. \label{eq:i_de_j}
\end{alignat}
In~\refQ,~$\lambda_i$ is an estimation of the expected computational cost of evaluating~\refP{\phi_i},~$p_{ij}$ is an estimation of the probability that~$c_j\leq0$ is satisfied at~$\phi_i$,~$r_{ij}$ is an estimation of the probability that~$\phi_i$ is representative for~$c_j\leq0$, and~$i:J\rightarrow I$ is a function that returns the lowest fidelity index~$i\in I$ for which~$\phi_i$ is representative for~$c_j\leq0$ for all feasible points in~$H$.
In~\cref{eq:objQ},~$f_Q(a)$ expresses an expected evaluation cost of~\refPbar, which depends on the evaluated fidelities.
As an illustrative example, suppose~$L=m=4$,~$a=[1,2,2,4]$ and~$\Phi(a)=[\phi_1,\phi_2,\phi_4]$.
Then,~$f_Q([1,2,2,4])=\lambda_1+\lambda_2p_{11}+\lambda_4p_{11}p_{22}p_{23}$.
To evaluate~$x\in H$, the evaluation at~$\phi_1$ of cost~$\lambda_1$ always occurs, then the evaluation at~$\phi_2$ of cost~$\lambda_2$ occurs if there was no interruption when evaluation at~$\phi_1$, i.e., if~$c_1(x,\phi_1)\leq0$, which happens with probability~$p_{11}$, and so on for~$\phi_4$.

Minimal expected evaluation costs are obtained by assigning constraints to low fidelities.
This is balanced with \cref{eq:rep_is_1} to ensure constraints are assigned to fidelities with maximal estimated representativity, therefore avoiding erroneous feasibility estimations.
Let us then study the effects of such errors on the optimization process.
\begin{itemize}
    \item \textbf{An infeasible point~$x\in X\backslash\Omega$ is deemed feasible.} This only occurs if~1 is not an element of~$\Phi(a)$, all constraints are satisfied at all fidelities in~$\Phi(a)$, and~$\bar{f}(x) \geq f^*$.
    Hence, an infeasible point can not become an incumbent solution.
    The only possible harm an infeasible point deemed feasible can cause occurs when a solver mistakenly uses it to compute search directions.
    \item \textbf{A feasible point~$x\in\Omega$ is deemed infeasible.} This can occur after any evaluation of~\refP{\phi} where $\phi<1$. The main issue this can cause is when~$x$ would, if correctly identified as feasible, become the new incumbent solution.
    Generally, misidentifying feasible points can cause the solver to omit them to compute search directions. This is particularly harmful if an omitted point has a good objective function value, missing an opportunity to point the solver towards a minimum.
\end{itemize}

The drawbacks are more important and more frequent in the second case compared to the first.
To mitigate this, \cref{eq:rep_is_1} imposes that a blackbox constraint can only be assigned to a fidelity that is representative for all feasible points in~$H$.

\subsection{Theoretical analysis of the assignment problem}
\label{sec:theorie_q}

Subproblem~\refQ is constructed to compute an assignment vector from a sample set~$H\subset X$ that ensures some results when~\refPbar evaluates~$x\in H$. These results are shown here.

\begin{lemma}
    \label{prop:r_croissant}
    When~$\text{H}\cap\Omega\neq\varnothing$, for each~$j\in J$,~$\{r_{ij}\}_{i\in I}$ is monotone increasing with respect to~$i$.
    \begin{proof}
        For each~$i^1,i^2\in I$ where~$i^1<i^2$, and for each~$j\in J$,
        \begin{alignat*}{2}
            & R_{ij} && :=\{x\in\text{H}\cap\Omega:\1(c_j(x,\phi_\ell)>0)=\1(c_j(x,1)>0)\,\forall\,\ell\in I,\,\phi_\ell\geq\phi_i\} \quad \forall\,i\in I \\
            \implies & r_{i^1j} && =\frac{\lvert R_{i^1j}\rvert}{\lvert\text{H}\cap\Omega\rvert} \qquad \text{ and } \qquad
            r_{i^2j}=\frac{\lvert R_{i^2j}\rvert}{\lvert\text{H}\cap\Omega\rvert}.
        \end{alignat*}

        Because~$\{\phi_i\}_{i\in I}$ is strictly increasing with respect to~$i$,~$\{\phi_\ell:\ell\in I,\phi_\ell\geq\phi_{i^2}\}\subset\{\phi_\ell:\ell\in I,\phi_i\geq\phi_{i^1}\}$, and therefore~$R_{i^1j}\subseteq R_{i^2j}$. As a result,~$\lvert R_{i^1j} \rvert\leq\lvert R_{i^2j} \rvert$ and~$r_{i^1j}\leq r_{i^2j}$.
    \end{proof}
\end{lemma}

\begin{theorem}
    \label{prop:cQ}
    For a given~$H\subset X$, when the fidelity controller performs an evaluation of~$x\in H\cap\Omega$ with an assignment vector given by Subproblem~\refQ, feasibility is conserved between~\refPbar and~\refP{1}.
    \begin{proof}
        If~$H\cap\Omega=\varnothing$, then~\refPbar is identical to~\refP{1}.
        Otherwise, \cref{eq:rep_is_1} ensures that each constraint~$c_j\leq0$ is assigned to a fidelity~$\phi_i$ where~$r_{ij}=1$, and Lemma~\ref{prop:r_croissant} indicates that~$r_{\ell j}=1$ for each~$\ell\geq i,\ell\in I$.
        A representativity of~1 means that the violation of a constraint is correctly identified at points in~$H\cap\Omega$.
        This implies that if a constraint is satisfied at the fidelity it is assigned to for a point~$x\in H^\text{ap}$, then it is also satisfied at any higher fidelity, including the truth.
    \end{proof}
\end{theorem}

Moreover, if~$x\in H\cap\Omega$ and~$1\in\Phi(a)$,~$\bar c(x) = c(x,1)$, and~$\bar f(x) = f(x,1)$.
\begin{assumption}
    \label{hyp:eval_supp}
    In Algorithm~\ref{algo:fico}, the evaluation of~\refP{1} that is conditional to~$1\notin \Phi(a) \text{ and } \bar f < f^*$ has no impact on the expected cost of evaluating~\refPbar.
\end{assumption}
This evaluation ensures that only feasible points can become incumbent solutions in the solver.
The assumption is false, but since the moments where a solution could become a new incumbent are unpredictable, it is necessary for the following theorem.

\begin{theorem}
    \label{prop:fQ}
    Under Assumption~\ref{hyp:eval_supp}, for a given~$H\subset X$ where~$H\cap\Omega\neq\varnothing$ and a given assignment vector~$a\in I^m$, the objective function~$f_Q(a)$ of Subproblem~\refQ expresses the expected cost of evaluating~\refPbar at a point~$x\in H^\text{ap}$.
    \begin{proof}
        For any~$H\subset X$ where~$H\cap\Omega\neq\varnothing$ and any~$a\in I^m$,
        \begin{alignat*}{2}
            \mathbb{E}_{x\in H^\text{ap}}\left[\text{cost of evaluation of~\refPbar at }x\right] = \sum_{i=1}^L\lambda_i & \mathbb{P}\left[\text{evaluation of~\refP{\phi_i} at }x\text{ occurs for an } x\in H^\text{ap}\right] && \\
             = \sum_{i:\phi_i\in\Phi(a)}\lambda_i\prod\limits_{\ell:\phi_\ell\in\Phi(a),\,\ell<i} & \mathbb{P}[\text{no interruption occurs at }(x,\phi_\ell)\text{ for an } x\in H^\text{ap}].
        \end{alignat*}
        This last result holds because the probability that an evaluation of~\refP{\phi} occurs is~0 if~$\phi\notin\Phi(a)$, and otherwise, it is the probability that no interruption happens beforehand.
        
        According to Theorem~\ref{prop:cQ}, considering the fidelity at which a constraint is assigned is sufficient to verify if this constraint would cause an interruption at higher fidelities.
        Therefore, assuming that all~$p_{a_jj}$ for each~$j\in J$ are independent and that the product of the elements of an empty set has value~1,
        \begin{alignat*}{1}
            \mathbb{P}[\text{no interruption occurs at }(x,\phi_\ell)\text{ for an } x\in H^\text{ap}] & = \prod_{j:a_j=\ell}p_{\ell j} \quad \forall \ell\in I \\
            \implies \mathbb{E}_{x\in H^\text{ap}}\left[\text{cost of evaluation of~\refPbar at }x\right] & = \sum_{i:\phi_i\in\Phi(a)}\lambda_i\prod\limits_{\ell:\phi_\ell\in\Phi(a),\,\ell<i}\left(\prod_{j:a_j=\ell}p_{\ell j}\right) \\
            & = \sum_{i:\phi_i\in\Phi(a)}\lambda_i\prod_{j:a_j<i}p_{a_jj} = f_Q(a).
        \end{alignat*}
    \end{proof}
\end{theorem}

A~priori constraints are ignored in Theorem~\ref{prop:fQ} because their violation causes the cost of an evaluation to be virtually null.
For a new candidate~$x\in X\backslash H$ that shares a similar behaviour with the points in~$H$,~$f_Q(a)$ is an approximation of the expected cost of~\refPbar, and this problem approximates the feasibility of~$x$.
As noted in~\cite{AlAuDiLedLe23}, richer sample sets~$H$ yield significantly better results.

\section{Two interruptible direct search algorithms}
\label{sec:inter_algos}

For the approximations based on~$H$ to be effective, information about many evaluated points at many fidelity levels must be available.
In this section, two optimization algorithms are proposed, each applicable under different assumptions regarding this availability.
The first is the Interruptible Direct Search (\ids) algorithm, which is presented in \cref{sec:ids}.
The second is its dynamic counter-part, the Dynamic Interruptible Direct Search (\dids) algorithm, which is presented in \cref{sec:dids}.
\dids relies on a stronger assumption, in order to exploit more fidelities and reach greater computational cost reductions than \ids.

\subsection{Interruptible Direct Search (\ids)}
\label{sec:ids}

The \ids algorithm is applicable under Assumption~\ref{hyp:LH}:

{\renewcommand{\theassumption}{2a}
\begin{assumption}
    \label{hyp:LH}
    A set of sample points~$\text{H}\subset X$, each evaluated at multiple fidelities with indices forming the set~$I$, is provided with Problem~\refP{\phi}.
    An evaluation of~\refP{\phi} only returns information for this fidelity level~$\phi$.
\end{assumption}}

A set~$H$ may be the result of Latin Hypercube Sampling (LHS) or previous experiments.
Using LHS ensures an unbiased sample set, but it may represent a high computational cost.
\ids is applicable to any problem illustrated in \cref{fig:classification} where \cref{hyp:LH} is verified.
The best solution from~$H$ is chosen as the initial optimization point.
Whenever the fidelity controller calls~\refP{\phi_i} at trial point~$x\in X$ for some~$i\in I$ during the optimization, it only receives the values~$f(x,\phi_i)$ and~$c(x,\phi_i)$ at cost~$\lambda_i$.
Fidelities with no assigned constraints are not used by \ids to avoid their cost.
Thus, for a given assignment~$a$, the set~$\Phi^\text{IDS}(a)$ is given by
\begin{alignat}{2}
    \Phi^\text{IDS}(a) = & \{\phi_i : i\in a\}. && \label{eq:Phi4.1}
\end{alignat}

The definition of~$\Phi^\text{IDS}(a)$ implies that for a single evaluation of~\refPbar, the fidelity controller could perform calls to~\refP{\phi} at several fidelities, with total cost exceeding the cost of evaluating~\refP{1}.
Conversely, using various fidelities may be relevant to interrupt evaluations on infeasible points.
The optimal solution of Subproblem~\refQ is the assignment vector~$a$ that finds the best compromise.

Although a differentiable formulation of Subproblem~\refQ exists~\cite{AlAuDiLedLe23}, it remains a mixed-integer problem with a polynomial objective function.
To help solving, Algorithm~\ref{algo:cuts} is proposed to add constraints to~\refQ that cut down the size of the set of feasible solutions, without removing an optimal solution.
Then, an exhaustive search on this feasible set is sufficient in practice to solve the optimal assignment problem.
This cutting algorithm introduces the set~$J_\nu\subseteq J$, given by~\eqref{eq:Jn}.
It is the set of constraint indices that are not significantly affected by multi-fidelity.
These constraints may have values that change with fidelity levels, but never such that their feasibility is affected.
\begin{equation}
    J_\nu = \{j\in J:\forall i\in I, p_{ij}=p_{L j} \text{ and } r_{ij}=1\}. \label{eq:Jn}
\end{equation}

\begin{algo}
\label{algo:cuts}
\begin{center}
    \begin{table}[ht]
    \centering \renewcommand{\arraystretch}{1.2}
    \begin{tabular}{L}
        \hline
        \textbf{Algorithm \thealgocounter: }\text{Cutting algorithm} \\
        \hline
        \textbf{Inputs: } \text{Subproblem~\refQ}; \text{subset of constraint indices }J_\nu \\
        \textbf{Outputs: } \text{Updated Subproblem~\refQ}\\
        \textit{1. Blackbox constraints related cuts} \\
        \left.\quad
            \begin{array}{l}
                \textbf{If } J=J_\nu,\,\textbf{Add }\text{constraints ``}a_j=1,\,\,\forall\,j\in J\text{'' to Subproblem~\refQ} \\
                \textbf{Else if }p_{L j}<1\text{ for some }j\in J_\nu \\
                \left.\middle|
                    \begin{array}{l}
                        \textbf{Add }\text{constraints ``}a_{\jc}\leq a_{j},\,\,\forall\,(\jc,j)\in J_\nu\times J \text{'' to Subproblem~\refQ} \\
                    \end{array}
                \right. \\
                \textbf{Else if }J_\nu\text{ is not empty} \\
                \left.\middle|
                    \begin{array}{l}
                        \textbf{Remove }\text{some decision variables of Subproblem~\refQ with }a\leftarrow\{a_j\in a:j\notin J_\nu\} \\
                    \end{array}
                \right. \\
            \end{array}
        \right. \\
        \textit{2. Fidelity related cuts} \\
        \left.\quad
            \begin{array}{l}
            \textbf{For each fidelity index } i \in I \\
            \left.\middle|
                \begin{array}{l}
                    \textbf{If }\text{there exists a fidelity index~$\ell\in I$,~$\ell>i$, such that~$\lambda_\ell\leq\lambda_i$} \\
                    \left.\middle|
                        \begin{array}{l}
                            \textbf{Add }\text{constraints ``}a_{j}\neq i,\,\,\forall\,j\in J \text{'' to Subproblem~\refQ} \\
                        \end{array}
                    \right. \\
                    \textbf{If } i \notin \bigcup\limits_{j\in J}i(j) \\
                    \left.\middle|
                        \begin{array}{l}
                            \textbf{Add }\text{constraints ``}a_{j}\neq i,\,\,\forall\,j\in J \text{'' to Subproblem~\refQ} \\
                        \end{array}
                    \right. \\
                \end{array}
            \right. \\
        \end{array}
    \right. \\
    \hline
    \end{tabular}
    \end{table}
\end{center}
\end{algo}

The first category of cuts from Algorithm~\ref{algo:cuts} follows from the observation that the feasibility of constraints of indices in~$J_\nu$ can be asserted at any fidelity~$\phi_i,i\in I$.
If all constraints are of this nature, they are simply assigned to the lowest fidelity,~$\phi_1$.
Else, if~$p_{L j}<1$ for some~$j\in J_\nu$, constraints of indices in~$J_\nu$ are certainly assigned to the lowest fidelity in~$\Phi^\text{IDS}(a)$.
Finally, if all constraints~$c_j\leq0,j\in J_\nu$ are such that~$p_{L j}=1$, they will never cause an interruption according to the parameters of Subproblem~\refQ.
As a result, in an optimal solution, these constraints can be assigned to any fidelity where a constraint~$c_j\leq0,\,j\notin J_\nu$ is also assigned, therefore not causing a new call to~\refP{\phi}.
Hence, Subproblem~\refQ is solved without assigning them, and they are assigned to the lowest fidelity where another constraint is assigned afterwards, as shown later in Algorithm~\ref{algo:IDS}.

The second category of cuts are related to fidelity levels. The first cuts have the effect of discarding fidelities that are not cheaper than higher fidelity levels. This causes the remaining~$\lambda_i$ values to be strictly increasing with respect to~$i \in I$.
Then, the last cuts follow from Theorem~\ref{prop:fidelity_filter}.

In the optimizations results presented in \cref{sec:results}, with the cutting algorithm, Subproblem~\refQ is solved in at most one second.

\begin{lemma}
    \label{prop:p_identiques}
    For any sample~$H\subset X$ and for each constraint index~$j\in J$,~$p_{ij}=p_{Lj}$ for each fidelity index~$i\in I$ where~$i\geq i(j)$.
    \begin{proof}
        For any sample~$H\subset X$ and for each constraint index~$j\in J$, \cref{eq:i_de_j} and Lemma~\ref{prop:r_croissant} imply~$r_{ij}=1$ for each fidelity index~$i\in I$ where~$i\geq i(j)$. For those indices, for each~$x\in H$,~$\1(c_j(x,\phi_i)>0)=\1(c_j(x,1)>0)$, meaning that the proportion of points where constraint~$c_j(x,\phi_i)$ $\leq0$ is violated is the same across all fidelity indices~$i\in I$ where~$i\geq i(j)$.
    \end{proof}
\end{lemma}

\begin{theorem}
    \label{prop:fidelity_filter}
    In the particular case of Subproblem~\refQ where~$\Phi^\text{IDS}(a)$ is given by~\eqref{eq:Phi4.1}, there exists an optimal solution such that all blackbox constraints are assigned to fidelity of indices in~$\bigcup\limits_{j\in J}i(j)$.
    \begin{proof}
        Two cases are considered.
        First, if there exists a constraint index~$j\in J$ such that~$p_{i(j)j}=0$, then constraint~$c_j\leq0$ is always violated for a feasible solution~$a\in I^m$.
        Consequently,~$a$ is optimal if this constraint is assigned to~$\phi_{i(j)}$, the lowest fidelity satisfying \cref{eq:rep_is_1}.
        
        In the second case, let~$a$ be a feasible solution where there exists a fidelity index~$\ic\in I\backslash\bigcup_{j\in J}i(j)$ such that at least one constraint is assigned to this fidelity.
        If~$a$ does not exist, then each feasible solution is such that all constraints are assigned to fidelity indices in~$\bigcup_{j\in J}i(j)$.
        Otherwise, let~$\hat{a}$ be a solution identical to~$a$, except that all constraints assigned to~$\phi_\ic$ are assigned to~$\phi_{\ic-1}$ instead, as defined by
        \begin{equation}
            \hat{a}_{j} =
            \begin{cases}
                \ic-1 \qquad & \text{ if }a_j=\ic \\
                a_{j} & \text{ otherwise}
            \end{cases} \qquad \forall\,j\in J. \label{eq:def_ac}
        \end{equation}
        Solution~$\hat{a}$ is feasible because for each~$j\in J$, if~$a_{j}=\ic$ and~$a$ is feasible, then~$\ic-1\geq i(j)$.
        Define~$P_i(a)=\prod_{j:a_j<i}p_{a_jj}$ for each~$i\in I$.
        Lemma~\ref{prop:p_identiques} implies that~$P_i(\hat{a})=P_i(a)$ for each~$i\in I\backslash\{\hat{i}\}$.
        Additionally, definition~\eqref{eq:def_ac} implies that~$\Phi^\text{IDS}(\hat{a})\backslash\{\ic-1,\ic\}=\Phi^\text{IDS}(a)\backslash\{\ic-1,\ic\}$, meaning that except for indices~$\ic-1$ and~$\ic$, the terms of the sum in the objective function~\eqref{eq:objQ} are identical for~$a$ and~$\hat{a}$. Moreover, if no constraint is assigned to~$\ic-1$ in~$a$,~$P_{\ic-1}(\hat{a})=P_\ic(a)$. Otherwise,~$P_{\ic-1}(\hat{a})=P_{\ic-1}(a)$ and the summation terms for~$i=\ic-1$ are also identical. As a result, with the indicator function~$\1(\ic-1\notin a)$ returning~1 if~$\ic-1\notin a$ and~0 otherwise,
        \begin{alignat*}{2}
            f_\mathcal{Q}(a)-f_\mathcal{Q}(\hat{a}) & = \lambda_\ic P_\ic(a)-\lambda_{\ic-1}P_{\ic-1}(\hat{a})\1(\ic-1\notin a) && = (\lambda_\ic-\lambda_{\ic-1}\1(\ic-1\notin a))P_\ic(a) \\
            & && > 0 \implies a\text{ is not optimal.}
        \end{alignat*}
        This inequality holds because~$P_\ic(a)=0$ only happens in the first case, and~$\lambda_i$ is strictly increasing with~$i\in I$ thanks to Algorithm~\ref{algo:cuts}. The contrapositive of this result is that a solution~$a$ is optimal if there exists no fidelity index~$\ic\in I\backslash\cup_{j\in J}i(j)$ such that at least one constraint is assigned to this fidelity.
    \end{proof}
\end{theorem}

If there is no constraint~$c_j\leq0,j\in J$, such that~$i(j)=L$, then the highest fidelity is not an element of~$\Phi^\text{IDS}(a)$, meaning that evaluations of the truth by the fidelity controller only happen for potential future incumbent solutions.
This allows IDS to save computational costs.
However, it may be preferable to ensure that~$1\in\Phi^\text{IDS}(a)$, such that points deemed feasible are always truly feasible.
Indeed, if~$1\in\Phi^\text{IDS}(a)$, the true values~$f(x,1)$ and~$c(x,1)$ are necessarily returned by~\refPbar when there is no interruption for point~$x\in X$.
To this effect, \ids posses the {\tt include\_truth} boolean parameter.
When it is true, the objective function~\eqref{eq:objQ} of the assignment Subproblem~\refQ is replaced with
\begin{equation}
    \lambda_L+\sum_{i:\phi_i\in\Phi^\text{IDS}(a),i\neq L}\lambda_i\prod_{j:a_j<i}p_{a_jj} \label{eq:force_eval_truth}
\end{equation}
to account for the mandatory evaluation using~$\phi_L=1$.
The \ids algorithm is shown in Algorithm~\ref{algo:IDS}.

\begin{algo}
\label{algo:IDS}
    \begin{center}
    \begin{table}[ht]
        \centering \renewcommand{\arraystretch}{1.2}
        \begin{tabular}{L}
            \hline
            \textbf{Algorithm \thealgocounter: }\text{Interruptible Direct Search (\ids)} \\
            \hline
            \textbf{Inputs: } \\
            \left.\middle|
               \begin{array}{rl}
                    \text{\refP{\phi}}: & \text{optimization problem defined by $f, c$ and $X \subseteq \R^n$;} \\
                    H \subset X: & \text{finite sample set;} \\
                    \textsf{solver}: & \text{direct search blackbox optimization solver;} \\
                    \text{\tt include\_truth}: & \text{boolean parameter that imposes~$1\in\Phi^\text{IDS}(a)$ when true. False by default.}
                \end{array}
            \right. \\
            \textbf{Output: }\text{ best solution found in }X \\
            \textbf{1. Optimal assignment vector computation} \\
            \left.\middle|
               \begin{array}{l}
                \text{Initialize $f^*\leftarrow\infty$, the incumbent value} \\
                \textbf{If }H\cap\Omega\text{ is empty, $a_j=L,\,\forall\,j\in J$ (equivalent to solving~\refP{1} without interruptions)} \\
                \textbf{Else } \\
                \left.\middle|
                   \begin{array}{l}
                        \text{From~$H$,~\eqref{eq:r},~\eqref{eq:p} and~\eqref{eq:kappa}, evaluate the parameters of Subproblem~\refQ}  \\
                        \textbf{If}\text{ {\tt include\_truth} is true, use~\eqref{eq:force_eval_truth} as objective function for
                        Subproblem~\refQ} \\
                        \text{Compute~$J_\nu$ using~\eqref{eq:Jn}, define Subproblem~\refQ and apply Algorithm~\ref{algo:cuts}} \\
                        \text{Find $a$ by performing an exhaustive search on the feasible set of Subproblem~\refQ} \\
                        \textbf{If }\text{constraints of index in $J_\nu$ are unassigned, }a_{j}\leftarrow\min_{j\in J\backslash{J_\nu}}\{a_j\}\text{ for each }j\in J_\nu \\
                    \end{array}
                \right. \\
                \end{array}
            \right. \\
            \textbf{2. Direct search optimization} \\
            \left.\middle|
               \begin{array}{l}
                    \text{Launch \textsf{solver} on ~\refPbar with the blackbox provided by Algorithm~\ref{algo:fico} (parametrized by $a$ and $f^*$)} \\
                    \quad\text{- Use the best point of $H$ as the initial point} \\
                    \quad\text{- Update $f^*$ after each new incumbent} \\
                    \quad\text{- If {\tt include\_truth} is true, use $\Phi^\text{IDS}(a)\cup\{1\}$ instead of $\Phi^\text{IDS}(a)$}\\
                \end{array}
            \right. \\
            \textbf{Return }\text{the \textsf{solver} output} \\
            \hline
        \end{tabular}
    \end{table}
    \end{center}
\end{algo}

As the incumbent value~$f^*$ varies during the optimization, Problem~\refPbar evaluated at the same point at different moments may return different values.
As a result,~$\bar f(x)$ and~$\bar c(x)$ require~$f^*$ as a third parameter to be functions, but they are considered as such to simplify the expressions.

\subsection{Dynamic Interruptible Direct Search (\dids)}
\label{sec:dids}

The \dids algorithm is applicable under an assumption regarding the accessibility of multi-fidelity information that differs from \ids.
The initial sample~$H$ is not necessary.
Assumption~\ref{hyp:LH} is replaced with Assumption~\ref{hyp:intermediaire}.
It is also later shown that Assumption~\ref{hyp:eval_supp} is no longer necessary.

{\renewcommand{\theassumption}{2b}
\begin{assumption}
    \label{hyp:intermediaire}
    Problem~\refP{\phi} has intermediary outputs for a set of multiple fidelity indices~$I$, as given by Definition~\ref{def:intermediaire}.
\end{assumption}}

\begin{figure}[htb]
    \begin{center}
    \begin{tikzpicture}[scale=1, transform shape]
    
        \node[extblock, minimum width=14.5cm, minimum height=2.75cm]                                                          (bb)     at (0,-1.75)   {};
        \node[draw, dashed, inner sep=0.3cm, rounded corners, minimum width=15.5cm, minimum height=8.75cm, anchor=north west] (bbw)    at (-0.5,6.7)  {};
        \node[extblock, minimum width=14.5cm, minimum height=4.2cm]                                                           (fico)   at (0,1.75)    {};
        \node[extblock, minimum width=10cm]                                                                                   (solver) at (0,7.33)    {\textbf{Solver}};
    
        \node[block]                      (forloop)     at ($(fico.west)+(1.5,-0.3)$)   {For each \\ $i\in I$ in \\ increasing \\ order};
        \node[decision, text width=2.8cm] (interrupter) at ($(forloop)+(11,0)$)         {};
        \node[decision, text width=1.4cm] (loopend)     at ($(forloop)+(2.75,0)$)       {};
        \node[block]                      (outputs)     at ($(fico.north)+(0.5,-1)$)    {$\bar f(x) = f(x,\phi)$ \\ $\bar c(x) = c(x,\phi)$};

        \node[decision, text width=1.4cm] (fid1)  at ($(bb.west)+(1.5,0.25)$)  {};
        \node[block]                      (eta1)  at ($(fid1)+(5,0)$)          {Pull $\eta_1$ samples of $f_{\xi}(x)$ and $c_{\xi}(x)$};
        \node[block]                      (etai)  at ($(eta1)+(0,-1)$)         {Pull $\eta_i-\eta_{i-1}$ samples of $f_{\xi}(x)$ and $c_{\xi}(x)$};
        \node[block]                      (saa)   at ($(eta1)+(6,0)$)          {Average $\eta_i$ samples \\ as in \cref{eq:saa}, \\ for $f$ and for $c$};
    
        \node[texte]          at ($(bbw.north)-(2,0.3)$)       {\textbf{Blackbox wrapper problem }\refPbar};
        \node[texte]          at ($(fico.north)+(0,-3.5)$ )    {\textbf{Fidelity controller algorithm}};
        \node[texte]          at ($(bb.north)+(0,-0.3)$ )      {\textbf{Blackbox problem }\refP{\phi}};
        \node[texte]          at ($(loopend)-(0,0.05)$)        {Is\\ $i=L$ \\ ?};
        \node[texte]          at ($(interrupter)-(0,0)$)       {Is \\ $c_j(x,\phi_i)>0$ \\ for a $j\in J$ where \\ $\phi_i\geq\phi_{a_j}$ \\ ?};
        \node[texte]          at ($(outputs)+(-2.4,0)$)        {$x$ is deemed \\ feasible};
        \node[texte]          at ($(outputs)+(2.8,0)$)         {interruption, \\ $x$ is deemed \\ infeasible};
        \node[texte]          at ($(fid1)-(0,0)$)              {Is \\ $i=1$ \\ ?};
        \node[texte] (inputs) at ($(solver.east)+(2,-0.33)$)   {$a$};
        
        \draw[arrow]  (loopend.west)          -- node[pos=0.4 ,above]  {no}     (forloop.east);
        \draw[arrow]  (interrupter.north)     |- node[pos=0.2 ,right]  {yes}    ($(outputs.east)+(0,0.25)$);
        \draw[arrow]  (interrupter.west)      -- node[pos=0.1 ,below]  {no}     (loopend.east);
        \draw[arrow]  (loopend.north)         |- node[pos=0.1 ,left]   {yes}    (outputs.west);

        \draw[arrow]  (fid1.east)    -- node[pos=0.4 ,above]         {yes}    (eta1.west);
        \draw[arrow]  (fid1.south)   |- node[pos=0.5 ,below]         {no}     (etai.west);
        \draw[arrow]  (eta1.east)    --                                       (saa.west);
        \draw[arrow]  (etai.east)    -|                                       (saa.south);
        
        \draw[arrow]          ($(solver.south)+(-4,0)$)   -- node[pos=0.1 ,left ] {$x$}                         ($(forloop.north)-(0.5,0)$);
        \draw[arrow]          ($(forloop.south)-(0,0)$)   -- node[pos=0.7 ,left ] {$x,\phi_i$}                  (fid1.north);
        \draw[arrow]          ($(outputs.north)+(1,0)$)   -- node[pos=0.8 ,right] {$\bar f(x),\bar c(x)$}       ($(solver.south)+(3.75,0)$);
        \draw[arrow]          (inputs.south)              --                                                    ($(fico.north)+(4.75,0)$);
        \draw[arrow]          (saa.north)                 -- node[pos=0.55 ,left ] {$ f(x,\phi_i),c(x,\phi_i)$}  (interrupter.south);
        \draw[arrow, dashed]  ($(solver.south)+(3.75,0)$) -- ($(solver.south)+(3.75,0.25)$) -- ($(solver.south)+(-4,0.25)$) -- ($(solver.south)+(-4,0)$);
        
    \end{tikzpicture}
    \end{center}
    \caption{Flow chart diagram showing a particular case of the fidelity controller with the \dids algorithm and a stochastic blackbox with SAA.}
    \label{fig:dynamic_fid_flow}
\end{figure}
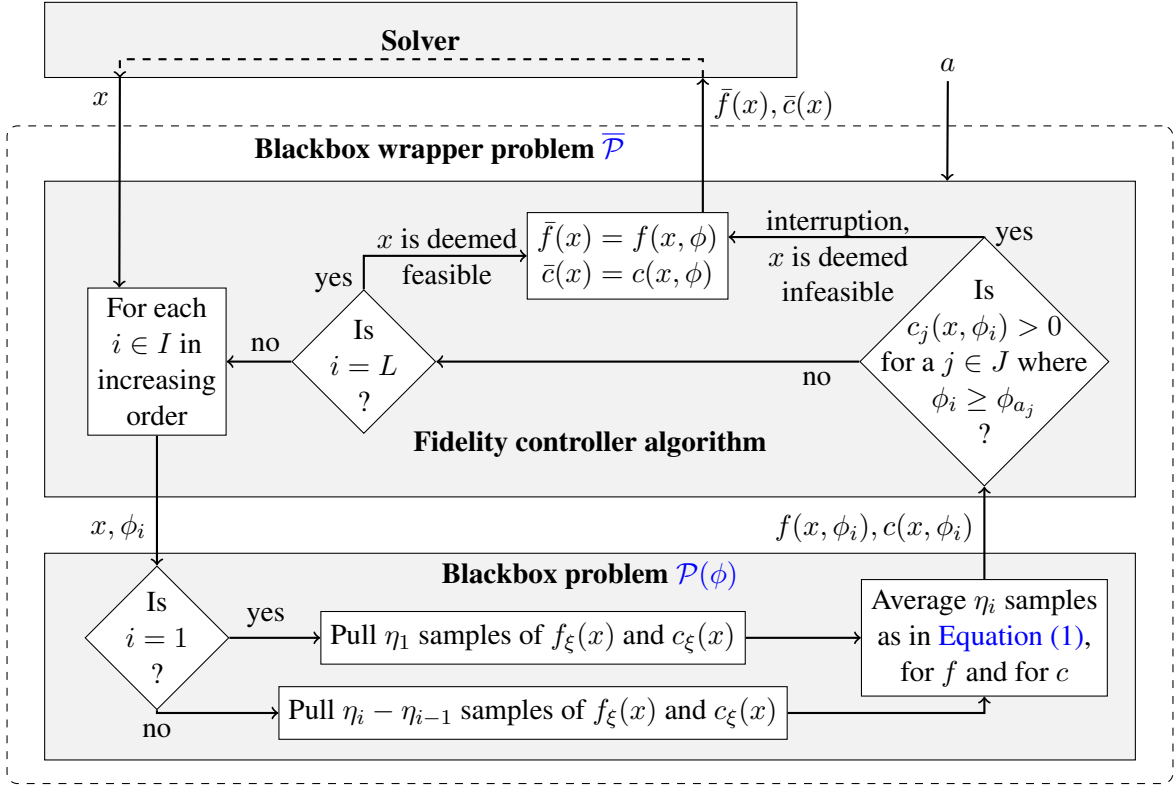

The fidelity controller algorithm is given a more specific interpretation for an evaluation using~$\phi_i,\,i\in I\backslash\{1\}$.
In the case of available intermediary outputs, an evaluation of~\refP{\phi_i} at trial point~$x\in X$ means that the blackbox continues its evaluation of~$x$ from~$\phi_{i-1}$ until fidelity~$\phi_i$ is reached.
\cref{fig:dynamic_fid_flow} illustrates this case where Definition~\ref{def:intermediaire} is met with SAA.
The cost~$\lambda_i$ represents the additional cost to reach~$\phi_i$ from~$\phi_{i-1}$.

Contrary to \ids, evaluating many fidelities cannot be more costly than evaluating~\refP{1}.
When evaluating~\refPbar at point~$x\in X$, fidelities are reached sequentially in increasing order, until~$\phi_L=1$ is reached.
Then,~$c(x,\phi_i)$ and~$f(x,\phi_i)$ are known for every fidelity~$\phi_i,\,i\in I$.
Consequently,
\begin{eqnarray}
    \Phi^\text{DIDS}(a) = \{\phi_i\}_{i\in I}
    & \quad \implies \quad & 
    f_\mathcal{Q}(a) = \sum_{i=1}^L\left(\lambda_i\prod_{j:a_j<i}p_{a_jj}\right).  \label{eq:f_dids}
\end{eqnarray}

This implies that a point is only deemed feasible when it is truly feasible.
As a result, Algorithm~\ref{algo:fico}'s safeguard against an infeasible becoming an incumbent is never run, and Assumption~\ref{hyp:eval_supp} is unnecessary.

Without Assumption~\ref{hyp:LH}, no initial sample~$H\subset X$ is required.
Rather, evaluated points deemed feasible during the optimization are stored in a cache~$\mathcal{V}$,
along with a vector of length~$m$ containing the lowest representative fidelity for each constraint for~$i(j)$ computations.
Then, periodically during the optimization process, points from~$\mathcal{V}$ are selected to form a sample set~$H$ from which~$a$ is updated.
Before iteration~$k$ of~\textsf{solver}, the smallest radius $\Delta$ such that~$B_\Delta(x^k)\cap\mathcal{V}\cap\Omega$ contains~$n+1$ feasible points is computed.
This ensures~$H$ is formed with points in close proximity to~$x^k$, while its size is at least the size of an~$n$-dimensional simplex.
When this radius does not exist, every constraint is assigned to~$\phi_L=1$.
The new optimal assignment~$a$ follows form the following theorem.

\begin{theorem}
    The particular case of Subproblem~\refQ where the objective function is given by~\eqref{eq:f_dids} has optimal solution~$a^*$ defined by
    \begin{equation}
        a^*_{j} = i(j). \label{eq:a_optimal}
    \end{equation}
    \begin{proof}
        Two cases are considered.
        First, if~$p_{i(j)j}=1$ for each~$j\in J$, then no blackbox constraint is ever violated for any feasible solution~$a\in I^m$. In that case, every feasible solution, notably~$a^*$, is optimal.
        If there exists a constraint index~$j\in J$ such that~$p_{i(j)j}=0$, then constraint~$c_j\leq0$ is always violated.
        Consequently, a solution~$a$ is optimal if this constraint is assigned to the lowest possible fidelity for a feasible solution, which is the case for~$a^*$.
        
        In the second case, let~$a$ be a feasible solution where there exists a constraint index~$\jc\in J$ such that~$p_{i(\jc),\jc}<1$ and~$a_{\jc}=\ic>i(\jc)$.
        If such an~$a$ does not exist, then~$a^*$ is the only feasible solution, and it is therefore optimal.
        Let~$\hat{a}$ be a solution identical to~$a$, except that constraint~$c_\jc\leq0$ is assigned to~$\phi_{i(\jc)}$ instead of~$\phi_\ic$, as defined by
        \begin{equation*}
            \hat{a}_{j} =
                \begin{cases}
                    i(\jc)\qquad & \text{ if }j=\jc \\
                    a_{j} & \text{ otherwise}
                \end{cases} \qquad \forall\,j\in J.
        \end{equation*}
        Solution~$\hat{a}$ is feasible because~$a$ satisfies
        \cref{eq:rep_is_1}.
        Define~$P_i(a)=\prod_{j:a_j<i}p_{a_jj}$ for each~$i\in I$.
        Definition~\eqref{eq:a_optimal} and Lemma~\ref{prop:p_identiques} imply that
        \begin{equation*}
            P_i(\hat{a}) =
                \begin{cases}
                    P_i(a)p_{i(\jc),\jc} & \text{ if } i\in\hat{I}:=\{i(\jc)+1,i(\jc)+2,\dots,\ic\} \\
                    P_i(a) & \text{ otherwise}
                \end{cases} \qquad \forall\,i\in I. 
        \end{equation*}

        As a result, except for indices in~$\hat{I}$, the terms of the sum in the objective function~\eqref{eq:f_dids} are identical for~$a$ and~$\hat{a}$.
        \begin{alignat*}{1}
            f_\mathcal{Q}(a)-f_\mathcal{Q}(\hat{a}) & = \sum_{i\in\hat{I}}\lambda_iP_i(a) - p_{i(\jc),\jc}\sum_{i\in\hat{I}}\lambda_iP_i(a) = (1-p_{i(\jc),\jc})\sum_{i\in\hat{I}}\lambda_iP_i(a)>0 \\
            & \implies a\text{ is not optimal}.
        \end{alignat*}
        The inequality holds because~$p_{i(j)j}>0$ for each~$j\in J$, implying that the summation is non zero, and~$p_{i(\jc),\jc}<1$.
        The contrapositive of this result is that a solution~$a$ is optimal if~$a_{j}=i(j)$ for each~$j\in J$ where~$p_{i(\jc),\jc}<1$, which describes~$a^*$.
    \end{proof}
\end{theorem}

As the assignment vector~$a$ varies during the optimization, Problem~\refPbar evaluated at the same point at different moments may return different values.
As a result,~$\bar f(x)$ and~$\bar c(x)$ require~$a$ as a third parameter to be functions, but they are considered as such to simplify the expressions.
The \dids algorithm is shown in Algorithm~\ref{algo:DIDS}.

\begin{algo}
\label{algo:DIDS}
    \begin{center}
    \begin{table}[ht]
        \centering \renewcommand{\arraystretch}{1.2}
        \begin{tabular}{L}
            \hline
            \textbf{Algorithm \thealgocounter: }\text{ Dynamic Interruptible Direct Search (\dids)} \\
            \hline
            \textbf{Inputs: } \\
            \left.\middle|
               \begin{array}{rl}
                    \text{\refP{\phi}}: & \text{problem containing $X \subseteq \R^n$, $f$ and $c$;} \\
                    x^0: & \text{optimization starting point (if \textsf{solver} requires it);} \\
                    \textsf{solver}: & \text{direct search blackbox optimization solver} 
                \end{array}
            \right. \\
            \textbf{Output: }\text{ best solution found in }X \\
            \textbf{1. Initialization} \\
            \left.\middle|
               \begin{array}{l}
                    \text{Initialize $\mathcal{V} \leftarrow \varnothing$, the cache of evaluated points} \\
                    \text{Initialize $k \leftarrow 1$, the \textsf{solver} iteration counter} \\
                    \text{Initialize $f^* \leftarrow \infty$, the incumbent value} \\
                    \text{Initialize $a\leftarrow a^0$, the initial assignment vector in which all constraints are assigned to~$\phi_L$} \\
                \end{array}
            \right. \\
            \textbf{2. Direct search optimization} \\
            \left.\middle|
               \begin{array}{l}
                 \text{Launch \textsf{solver} on~\refPbar with the blackbox provided by Algorithm~\ref{algo:fico} (parametrized by $a$ and $f^*$)} \\
                \quad\text{- Use $x^0$ as the initial point} \\
                \quad\text{- After each evaluation of~\refP{\phi}, update $\mathcal{V}$} \\
                \textbf{After }\text{iteration~$k$ of \textsf{solver} with incumbent solution $x^k$} \\
                \left.\middle|
               \begin{array}{l}
                    f^*\leftarrow f(x^k) \\
                    \text{Find the minimal radius $\Delta$ such that the ball $B_\Delta(x^{k})\cap\mathcal{V}\cap\Omega$ contains at least $n+1$ points.} \\
                    \textbf{If }\text{$\Delta$ does not exist, } a\leftarrow a^0 \\
                    \textbf{Else} \\
                    \left.\middle|
                        \begin{array}{l}
                            H\leftarrow B_\Delta(x^{k})\cap\mathcal{V}  \\
                            \text{From~$H$, compute the representativity probability~$r_{ij}$ for each~$(i,j)\in I\times J$ with~\eqref{eq:r}} \\
                            a\leftarrow [i(j)]_{j\in J} \text{, with~\eqref{eq:i_de_j}, as indicated by~\eqref{eq:a_optimal}}
                        \end{array}
                    \right. \\
                    k \leftarrow k+1 \\
                \end{array}
            \right. \\
            \end{array}
            \right. \\
            \textbf{Return }\text{the \textsf{solver} output} \\
            \hline
        \end{tabular}
    \end{table}
    \end{center}
\end{algo}

The set~$H$ is updated dynamically such that it contains points in close proximity to~$x^k$, while~$H\cap\Omega$ is at least the size of a simplex.

\section{Computational results}
\label{sec:results}

In this section, the proposed  algorithms are compared on blackbox optimization benchmark problems.
Since the blackbox simulator from the Hydro-Québec project motivating this study is still under development, benchmarks are conducted on the \solar{} collection of test problems~\cite{solar_paper}.
To the authors' knowledge, the~\solar{} family of benchmarks is the only one containing realistic industrial problems in which a fidelity parameter affects the constraint values.
Moreover, selecting benchmark problems that reflect real engineering optimization challenges is crucial. 
This suite provides ten blackbox simulators of a concentrated solar power plant (CSP).
\cref{tab:solar} shows relevant characteristics for the three problems selected for this study:~\solar{2},~\solar{3} and~\solar{4}.
Problem~\solar{7} is the only other problem with multi-fidelity constraints.
However, preliminary tests show that the optimal assignment vector computed by \ids and \dids is almost always~$a_j=L$ for each~$j\in J$.
This means that fidelities that allow relevant interruptions are rarely found, and both algorithms almost always simply solve~\refP{1} directly.
Problem~\solar{7} is not studied further.
\begin{table}[ht!]
    \renewcommand{\tabcolsep}{4.5pt}
    \renewcommand{\arraystretch}{1.2}
    \begin{footnotesize}
    \begin{center}
        \begin{tabular}{|l|cc|c|cc|c|}
        \hline
        Problem & \multicolumn{3}{c|}{Number of variables} & \multicolumn{3}{c|}{Number of constraints} \\
         & continuous & integer (categorical) & $n$ & a~priori & simulated (multi-fidelity) & $m$ \\
        \hline
        \hline
        \solar{2} & 12 & 2\quad (0) & 14 & 5 & 7 \quad (4) & 12 \\
        \solar{3} & 17 & 3\quad (1) & 20 & 5 & 8 \quad (5) & 13 \\
        \solar{4} & 22 & 7\quad (1) & 29 & 7 & 9 \quad (6) & 16 \\
        \hline
        \end{tabular}
    \end{center}
    \end{footnotesize}
    \caption{\footnotesize Number of variables and constraints for the three studied \solar{} problems. In parentheses, the number of categorical variables and multi-fidelity constraints are indicated.}
    \label{tab:solar}
\end{table}

To conduct extensive benchmarks, the computational cost of using SAA with the selected \solar{} problems showed to be prohibitive.
Instead, the \solar{} problems are used deterministically by setting the seed value to 0.
For \dids, intermediary outputs are made available by adding code to \solar{}.
To this effect, the fidelity parameter, described in~\cite{solar_paper}, controls the convergence criteria of numerical methods.
In PRIAD, evaluation times may be measured in days~\cite{PRIAD_KoMeCoGaVoAlDeBl2021}.
To emulate this context where solver computing times are negligible, only blackbox evaluation times are considered.
Time dependant data profiles are used~\cite{MoWi2009,G-2025-36,Beiranvand2017} for comparisons on numerous problem instances.
For an algorithm, a data profile plots the portion of~$\tau$-solved problem instances with respect to time.
A problem instance is said to be~$\tau$-solved at time~$T$ by algorithm~$A$ with initial point~$x^0$ if
\begin{equation*}
    \frac{f(x^0)-f^T}{f(x^0)-f^*} \geq 1 - \tau,
\end{equation*}
where~$f^T$ is the best feasible objective function value found by algorithm~$A$ as of time~$T$,~$f^*$ is the best feasible solution found among all algorithms on the problem instance, and~$\tau\in(0,1)$ is a tolerance.

The interruption algorithms are paired with the \nomad blackbox optimization solver, version~4.4.0~\cite{nomad4paper}.
\nomad is a freely available implementation of the \mads algorithm~\cite{AuDe2006,AuDe09a} that has shown to be successful on real engineering and industrial problems~\cite{AlAuGhKoLed2020}.
It is among the best suited solvers for constrained, noisy and highly discontinuous problems where direct-search methods are preferred to model-based approaches~\cite{G-2025-70}.
The \nomad parameters are left to their default values for all tested algorithms, 
    except that the extreme barrier is used for all constraints to ensure the rejection of points for which the evaluation is interrupted by the fidelity controller.
The fidelity indices used in this section are~$I=\{1,2,\dots,11\}$ with corresponding fidelity parameters~$\{\phi_i\}_{i\in I}=\{10^{-10},0.1, 0.2, 0.3, \dots, 1\}$.
A unique feasible initial solution is provided for each \solar{} problem.

\ids and \dids are compared to two methods: a base case where \nomad solves~\refP{1} directly without the fidelity controller, and the \inter{DS} algorithm from~\cite{AlAuDiLedLe23}, also paired with \nomad.
For \ids and \inter{DS}, the set~$H$ is obtained via Latin Hypercube Sampling (LHS) with 1000 points for each problem, and the LHS bounds are selected close to the starting point, as recommended in~\cite{AlAuDiLedLe23}.
For \dids, the LHS set is not considered.
Rather, intermediary outputs are available, as per Definition~\ref{def:intermediaire}.
As the LHS is included in the \inter{DS} algorithm and \ids assumes the LHS results already exist, no LHS computation time is considered for fair comparison.

First, \cref{sec:s3_results} shows benchmarks on \solar{3}.
Second, \cref{sec:s4_results} presents results for \solar{4}, a problem where constraint behaviour with respect to fidelity is highly correlated with the blackbox input.
Last, \cref{sec:s2_results} studies benchmarks on \solar{2}, a problem where the impact of erroneous feasibility estimations by the fidelity controller defies expectations.
All optimizations are performed using multiple Intel Xeon Gold~6150 CPUs operating at 2.70~GHz.

\subsection[SOLAR3 benchmarks]{\solar{3} benchmarks}
\label{sec:s3_results}

\begin{figure}[ht!]
    \centering
    \begin{subfigure}[t]{0.49\textwidth}
        \includegraphics[width=\textwidth]{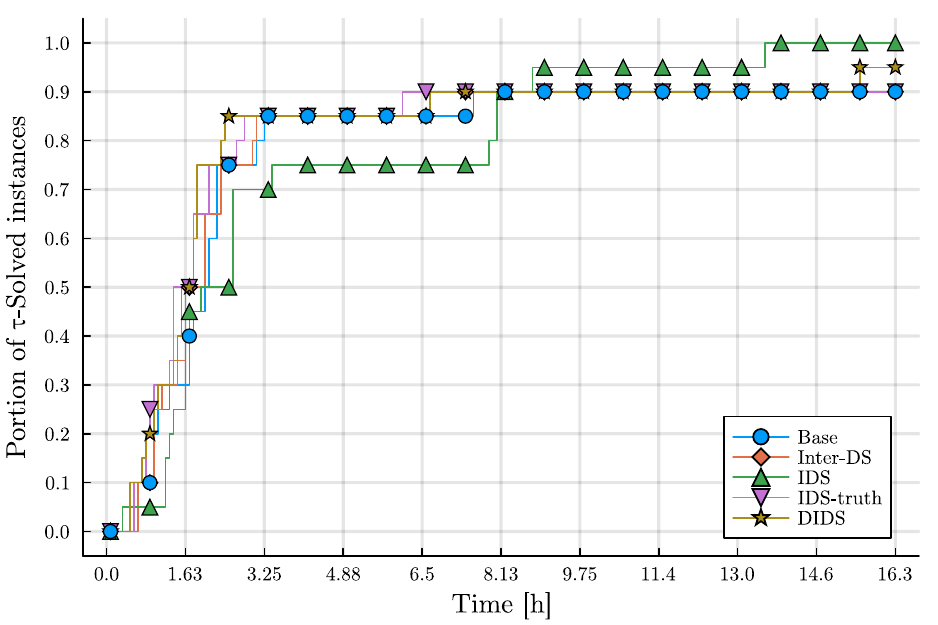}
        \caption{$\tau=0.2$}
        \label{fig:s3_0.2}
    \end{subfigure}
    \hfill
    \begin{subfigure}[t]{0.49\textwidth}
        \includegraphics[width=\textwidth]{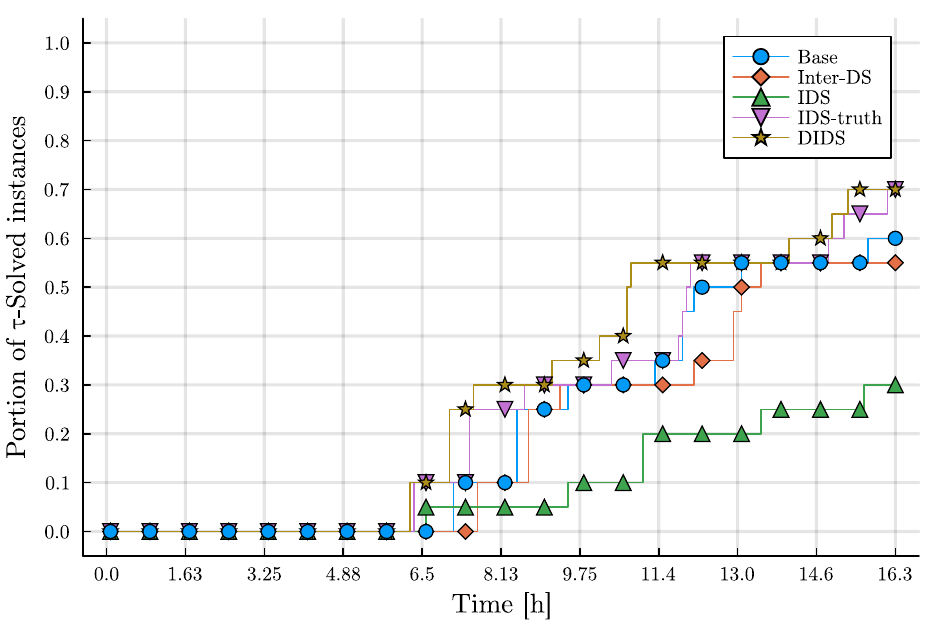}
        \caption{$\tau=0.01$}
        \label{fig:s3_0.01}
    \end{subfigure}
    \caption{Data profiles on Problem \solar{3} with 20 different problem instances.}
    \label{fig:dp_s3}
\end{figure}

The optimization results of 20 \solar{3} problem instances, obtained with different \nomad seeds and a $16.3$ hours budget, are shown in \cref{fig:dp_s3} with two different values of $\tau$.
A first observation is that pairing \nomad with \dids or \ids with {\tt include\_truth} set to true necessarily improves its performance.
At any time, with $\tau=0.01$, \dids~$\tau$-solves the most problem instances and \ids-truth~$\tau$-solves the second most. 
With tolerance $\tau=0.2$, \ids performs the best with $100\%$ of instances $\tau$-solved.
With~$\tau=0.01$, this percentage drops to $30\%$, which is significantly lower than the $70\%$ of both \dids and \ids-truth. 
\inter{DS} presents the same data profile as the base case with $\tau=0.2$, and $\tau$-solves $5\%$ less problem instances than the base case with the lower tolerance.

\begin{figure}[ht!]
    \centering
    \includegraphics[width=0.7\textwidth]{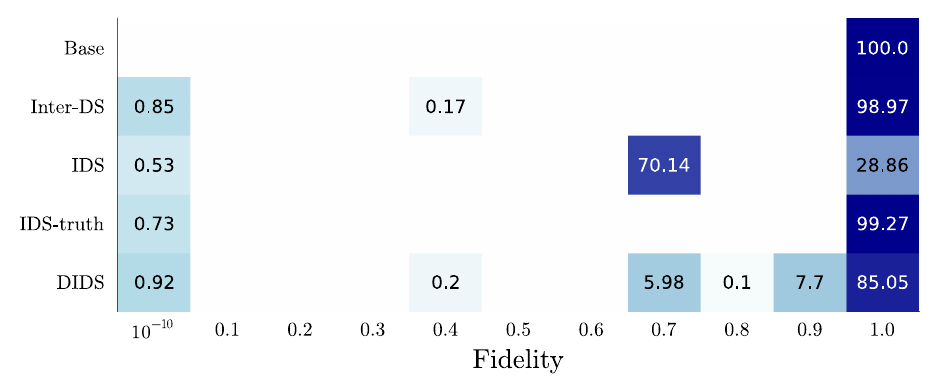}
    \caption{Occurrence (\%) of last used fidelities in \solar{3} evaluations.}
    \label{fig:s3phiprop_EB}
\end{figure}

To delve further into these results, \cref{fig:s3phiprop_EB} illustrates the frequency at which evaluations throughout all 20 optimizations ended.
In \solar{3}, the constraint that has lead to the most interruptions for all algorithms throughout all optimizations is~$c_2\leq0$, which verifies compliance to an energy demand.
For this constraint, higher fidelities best indicate if the energy produced is sufficient.
As a result, low fidelities are rarely used.
\dids leverages six different fidelity levels, resulting in the best performance with~$\tau=0.01$.
\ids only uses~$\phi_L=1$ in~$28.86\%$ of its evaluations, while the other algorithms use it for at least~$85.05\%$ of their evaluations.

As discussed at the end of \cref{sec:compute_a}, \ids and \dids are designed to ensure low probabilities that feasibility errors occur.
Ideally, an interruption method performs more evaluations than the base case in the same time budget, and interruptions always occur on infeasible points.
However, this can not be guaranteed without additional hypotheses, and despite low probabilities of errors, they are expected to occur during the optimization process (16.3 hours).
Such errors sometimes cause the evaluation sequence to differ from the base case.
To study this phenomenon, \cref{tab:info_add_s3} displays, for each interruption algorithm, the average factor of number of base case evaluations, the percentage of problem instances where the sequence of points up to the base case's last evaluation are different, and, for these instances, the average percentage of base case evaluations before the first different point in the sequence.
\begin{table}[ht!]
    \begin{footnotesize}
    \begin{center}
        \begin{tabular}{|l|cccc|}
        \hline
        Value compared with the base case & \inter{DS} & \ids & \ids-truth & \dids \\
        \hline
        \hline
        Average number of evaluations factor & 1.001 & 1.006 & 1.003 & 1.013 \\
        \% of problem instances where an $x$ differs & 0 & 100 & 0 & 0 \\
        Average \% of evaluations before an $x$ differs & - & 2.26 & - & - \\
        \hline
        \end{tabular}
    \end{center}
    \end{footnotesize}
    \caption{Evaluation sequence comparison with the base case for \solar{3}.}
    \label{tab:info_add_s3}
\end{table}

Where all other algorithms followed the same evaluation sequence as the base case for all optimizations, \ids has differences for all~$20$ optimizations.
These differences started after only~$2.26\%$ of the base case evaluations on average.
This is the result of \ids's prevalent usage of fidelities inferior to~$\phi_L=1$.
This way, \ids performs the most average number of evaluations at the cost of causing numerous feasibility errors, and these errors lead to worse search directions, resulting in the worst performance with~$\tau=0.01$.
Additionally, \ids's average number evaluations is only~$1.006$ times that of the base case.
Activating the {\tt include\_truth} parameter for \ids successfully prevents these feasibility errors. 
With~$\tau=0.01$, \inter{DS}'s profile is below the base case's at any time.
This is because~$98.97\%$ of its evaluation cost~$\lambda_1+\lambda_5+\lambda_{11}$.
The use of fidelity~$\phi_5=0.4$ is not sufficiently effective at causing interruptions to overcome the~$\lambda_5$ cost, causing a constant delay in the optimization.

\subsection[SOLAR4 benchmarks]{\solar{4} benchmarks}
\label{sec:s4_results}

\begin{figure}[ht!]
    \centering
    \begin{subfigure}[t]{0.49\textwidth}
        \includegraphics[width=\textwidth]{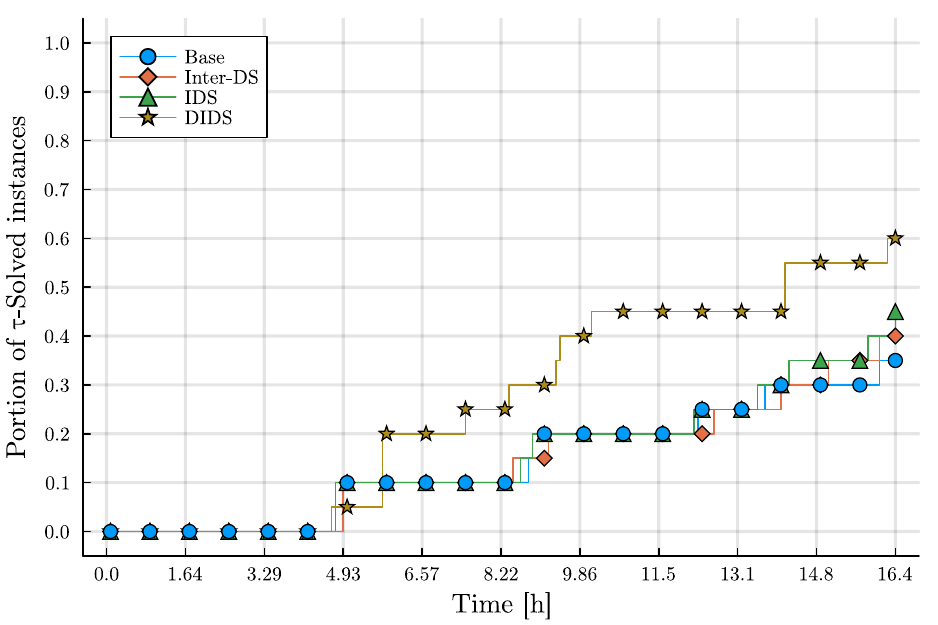}
        \caption{$\tau=10^{-3}$}
        \label{fig:s4_10-3}
    \end{subfigure}
    \hfill
    \begin{subfigure}[t]{0.49\textwidth}
        \includegraphics[width=\textwidth]{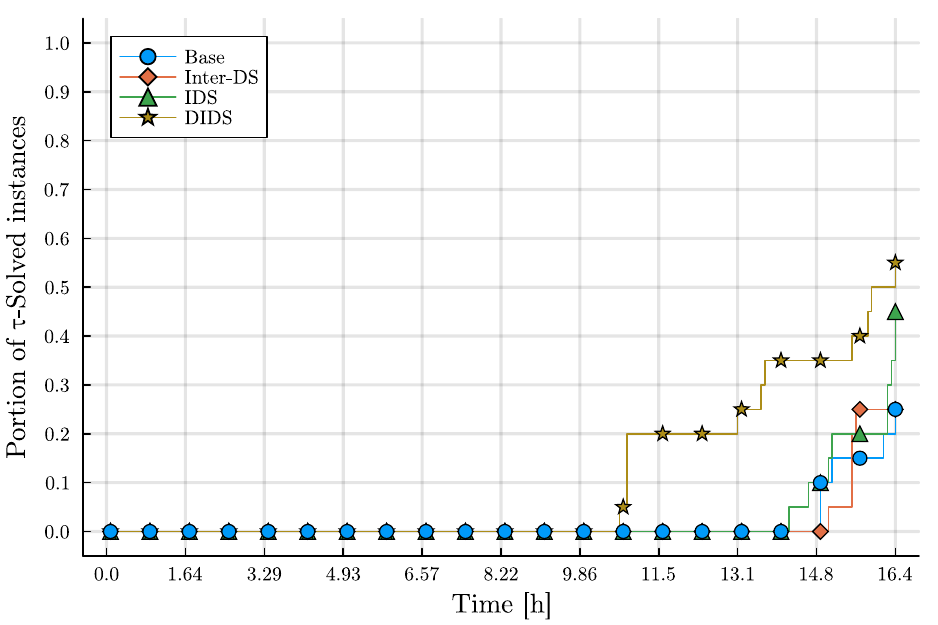}
        \caption{$\tau=10^{-5}$}
        \label{fig:s4_10-5}
    \end{subfigure}
    \caption{Data profiles on Problem \solar{4} with 20 different problem instances.}
    \label{fig:dp_s4}
\end{figure}

The optimization results of 20 \solar{4} problem instances, obtained with different \nomad seeds and a $16.4$ hours budget, are shown in \cref{fig:dp_s4} with two different values of $\tau$.
For \solar{4}, the assignment vector found by \ids contains~$\phi_L=1$, meaning that the {\tt include\_truth} parameter is redundant.
At any time, the data profiles of \inter{DS}, \ids and the base case are almost identical, even with low tolerance values such as $10^{-3}$ and $10^{-5}$ (which usually help differentiating the profiles).
The only exceptions to this similitude are the last few minutes with tolerance $\tau=10^{-5}$, where \ids suddenly reached $45\%$ of $\tau$-solved problem instances, compared to $25\%$ for \inter{DS} and the base case.
Conversely, \dids dominates both plots.

\begin{figure}[ht!]
    \centering
    \includegraphics[width=0.7\textwidth]{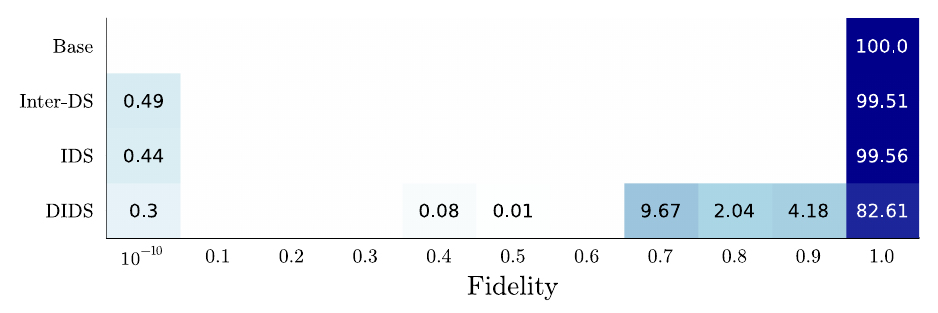}
    \caption{Occurrence (\%) of last used fidelities in \solar{4} evaluations.}
    \label{fig:s4phiprop_EB}
\end{figure}

Further insights are found in \cref{fig:s4phiprop_EB}.
Similarly to \solar{3}, the constraint that has lead to the most interruptions for all algorithms throughout all optimizations is~$c_2\leq0$, which verifies compliance to an energy demand and requires high fidelities to estimate.
\inter{DS} and \ids use~$\phi_L=1$ for~$99.51\%$ and~$99.56\%$ of their evaluations respectively.
In that matter, they behave almost identically to the base case.
The leveraging of intermediary outputs by \dids allowed it to interrupt evaluations at various fidelities. 

\cref{tab:info_add_s4} describes some evaluation sequence statistics for these optimizations.
Notice that \inter{DS} achieves less evaluations than unpaired \nomad on average.
Overall, the ineffectiveness of the non-dynamic interruption algorithms at levering multi-fidelity information is attributed to the fact that the behaviour of the constraints relative to fidelity is heavily dependant on the solution space for \solar{4}~\cite{AlAuDiLedLe23}.
As such, the assignment vectors of \inter{DS} and \ids can become ineffective, whereas the dynamic adaptations of \dids ensure a high quality assignment vector during the whole optimization.
This explains why \dids is the only algorithm that leverages multiple fidelity levels other than~$10^{-10}$ and~$1$ in \cref{fig:s4phiprop_EB}.
This use of multi-fidelity for interruptions causes \dids to not follow the base case's evaluation sequence for~$90\%$ of the problem instances, and the first difference occurs after only~$6.61\%$ of base case evaluations on average.
Nonetheless, deeper analysis shows the different search directions result in similar solutions to those of the base case.
Considering all factors, the reason for \dids's highest portion of~$\tau$-solved instances is the~$1.211$ average factor in number of evaluations.
\begin{table}[ht!]
    \begin{footnotesize}
    \begin{center}
        \begin{tabular}{|l|ccc|}
        \hline
        Value compared with the base case & \inter{DS} & \ids & \dids \\
        \hline
        \hline
        Average number of evaluations factor & 0.991 & 1.039 & 1.211 \\
        \% of problem instances where an $x$ differs & 70 & 70 & 90 \\
        Average \% of evaluations before an $x$ differs & 10.91 & 10.37 & 6.61 \\
        \hline
        \end{tabular}
    \end{center}
    \end{footnotesize}
    \caption{Evaluation sequence comparison with the base case for \solar{4}.}
    \label{tab:info_add_s4}
\end{table}

\subsection[SOLAR2 benchmarks]{\solar{2} benchmarks}
\label{sec:s2_results}

\begin{figure}[ht!]
    \centering
    \begin{subfigure}[t]{0.49\textwidth}
        \includegraphics[width=\textwidth]{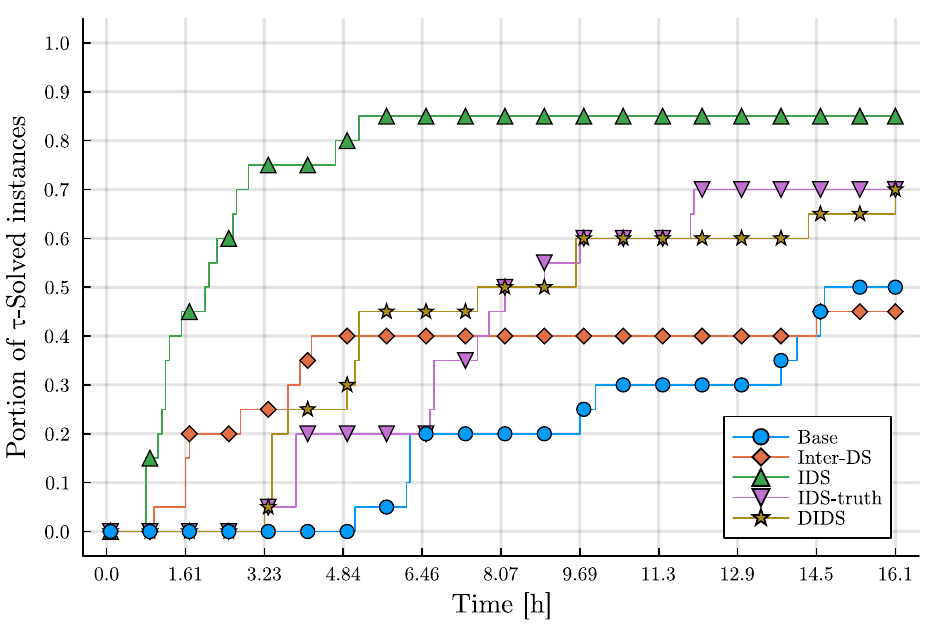}
        \caption{$\tau=0.1$}
        \label{fig:s2_0.1}
    \end{subfigure}
    \hfill
    \begin{subfigure}[t]{0.49\textwidth}
        \includegraphics[width=\textwidth]{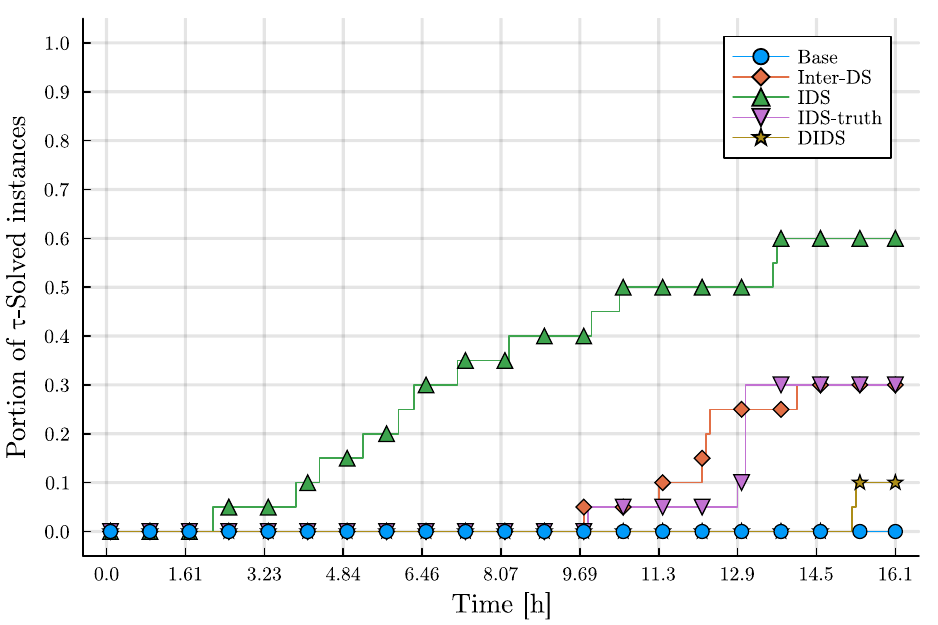}
        \caption{$\tau=0.01$}
        \label{fig:s2_0.01}
    \end{subfigure}
    \caption{Data profiles on Problem \solar{2} with 20 different problem instances.}
    \label{fig:dp_s2}
\end{figure}

The optimization results of 20 \solar{2} problem instances, obtained with different \nomad seeds and a $16.1$ hours budget, are shown in \cref{fig:dp_s2} with two different values of~$\tau$.
All interruption methods except \inter{DS} allow \nomad to reach better or equal solutions than by directly solving~\refP{1} at any time.
The \ids algorithm dominates both plots by $\tau$-solving $85\%$ of instances with a tolerance of $0.1$ and $60\%$ of instances with a tolerance of $0.01$.
By decreasing the tolerance, \dids and the base case are the algorithms with the largest decrease of $\tau$-solved instances, going from $70\%$ to $10\%$ for \dids and from $50\%$ to $0\%$ for the base case.

\begin{figure}[ht!]
    \centering
    \includegraphics[width=0.7\textwidth]{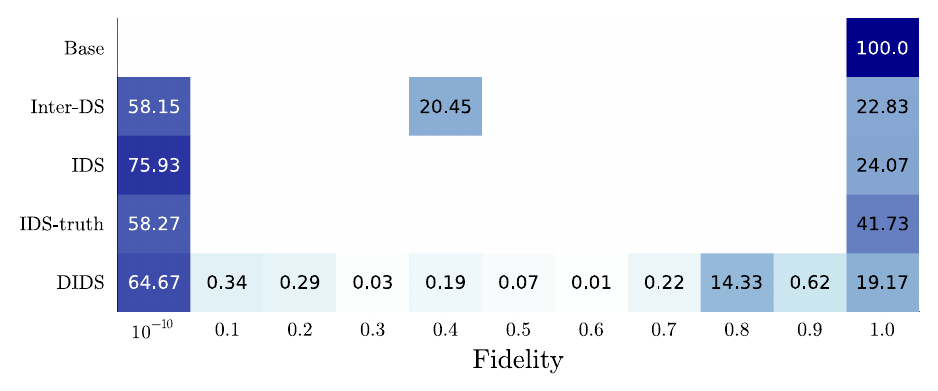}
    \caption{Occurrence (\%) of last used fidelities in \solar{2} evaluations.}
    \label{fig:s2phiprop_EB}
\end{figure}

To explain these results, \cref{fig:s2phiprop_EB} illustrates 
how each algorithm used the available fidelities.
It indicates that \ids interrupted evaluations at the lowest fidelity the most: $75.93\%$ of evaluations interrupted at~$\phi_1=10^{-10}$.
Throughout all problem instances and all algorithms, every interruption at~$\phi_1=10^{-10}$ was caused by~$c_6\leq0$, and by it only.
It ensures that the number of heliostats to place,~$x_6$, can fit in the field, for which the dimensions are computed from~$x_1$,~$x_2$ and~$x_3$.
This is not an a~priori constraint, but its value is computed before any ray-tracing MC draw is realized.
This explains why~$\phi_1=10^{-10}$ is so heavily used by all interruption algorithms.
The usage of all fidelity levels by \dids is not sufficient to overcome the other interruption methods.
This is not surprising, as the behaviour of the constraints relative to fidelity is fairly constant throughout the solution space of \solar{2}~\cite{AlAuDiLedLe23}, meaning that dynamic adjustments to the assignment vector are almost irrelevant.

\cref{tab:info_add_s2} displays how the different multi-fidelity approaches affected the optimizations.
A direct correlation between the frequency of evaluations interrupted at~$\phi_1=10^{-10}$ in \cref{fig:s2phiprop_EB} and the number of evaluations in this table is observed.
This is the first element explaining \ids's success: it performed~6.86 times more evaluations than the base case.
The second is that \ids's evaluation sequence differed from the base case for all problem instances, and it does so the fastest: only the first~$3.84\%$ of base case evaluations are identical before the first difference on average.
A deeper analysis revealed that the feasibility mistakes resulted in \nomad exploring better solutions.
This challenges the idea that interruption methods should be designed to minimize feasibility estimation errors.
Most errors occurred on infeasible points~$x\in X$ where~$f(x,1)$ is close to~$f^*$ and the true~$h(x)$ value is close to~$0$.
The fact that this lead to better solutions is most likely the result of chance.
Data profiles with a greater number of problem instances are required to confirm this observation.
\begin{table}[ht!]
    \begin{footnotesize}
    \begin{center}
        \begin{tabular}{|l|cccc|}
        \hline
        Value compared with the base case & \inter{DS} & \ids & \ids-truth & \dids \\
        \hline
        \hline
        Average number of evaluations factor & 4.05 & 6.86 & 2.94 & 4.65 \\
        \% of problem instances where an $x$ differs & 100 & 100 & 0 & 5 \\
        Average \% of evaluations before an $x$ differs & 5 & 3.84 & - & 29 \\
        \hline
        \end{tabular}
    \end{center}
    \end{footnotesize}
    \caption{Evaluation sequence comparison with the base case for \solar{2}.}
    \label{tab:info_add_s2}
\end{table}

In contrast, \dids differed from the base case in only~$5\%$ of instances, and for this~$5\%$, the first~$29\%$ of base case evaluations are identical before the first difference on average.
It also performed~$4.65$ times more evaluations than the base case on average.
Yet, \dids $\tau$-solves only~$10\%$ of the problem instances with tolerance~$\tau=0.01$.
Again, more tests would be needed to confirm these results.
Activating the {\tt include\_truth} parameter in \ids decreases the number of optimizations with a different sequence from~$100\%$ to~$0\%$, all while performing~$2.94$ times more evaluations than the base case on average.
Similarly to \solar{3}, this demonstrates that the parameter is efficient at reducing feasibility errors, but contrary to \solar{3}, \ids performs better here without the parameter.
This is because it reduces the number of evaluations factor from~$6.86$ to~$2.94$, but \nomad's performance is still improved when using \ids with the parameter activated.
In comparison, \inter{DS} performed more evaluations, with a factor of~$4.05$, but also had numerous feasibility errors (all problem instances differed from the base case), resulting in a similar data profile with tolerance~$\tau=0.01$.

\subsection{Remarks on the computational effort}
Recall that solver computing times are omitted in the data profiles to create conditions similar to PRIAD's.
This explains the irregular optimization times in \cref{fig:dp_s3,fig:dp_s4,fig:dp_s2}, as they are computed a~posteriori.
Note that a total of 6333 hours (8.65 months, divided between multiple machines) of optimization time was required to generate these figures.
To this effect, only pairings with \nomad are tested.
Numerical tests in~\cite{AlAuDiLedLe23} show that solvers implementing the Particle Swarm Optimization and the Nelder-Mead algorithms are successfully used with fidelity controller-based algorithms.

When \inter{DS} was introduced in~\cite{AlAuDiLedLe23}, it was noted that when~$H$ contains no feasible points, \inter{DS} is extremely harmful to the optimization.
These pathological cases are solved in \ids and \dids by assigning all constraints to~$\phi_L=1$ when~$H\cap\Omega=\varnothing$.
Because the solution is trivial, no additional test is conducted.

\section{Discussion}

This work introduces two new algorithms for constrained multi-fidelity blackbox optimization problems based on interruptible evaluations.
They both perform feasibility estimations at various fidelity levels, and interrupt the process when a point is deemed infeasible.
The first, \ids, exploits a set of points evaluated prior to the optimization to construct its interruption mechanism.
The second, \dids, relies on intermediary information from previously evaluated points during the optimization.
\ids is presented as a theoretical improvement of \inter{DS}, which is also based on pre-optimization evaluations.
The results show that for most problems, \ids also successfully improves on \inter{DS} practically.

For the problems tested, \ids with the {\tt include\_truth} parameter activated and \dids both strictly allow \nomad to reach better solutions at the end of the same time budget.
Without the {\tt include\_truth} parameter, \ids's performance is highly variable because it performs interruptions more aggressively.
The \solar{3} Problem punishes this greediness, \solar{4} is indifferent to it, and  \solar{2} rewards it.
A shortcoming of this study is that it is still unclear whetter feasibility estimation mistakes are desirable for \solar{2}, or if more optimizations should compose the data profiles to expose different global trends.
Future work will address phenomenon, as well as study if a bias is introduced by successful interruptions at low fidelities, leading to a lack of high fidelity information for infeasible points, which is used for assignment vector updates.

An important takeaway is that the behaviour of the constraints relative to fidelity must be fairly constant throughout the solution space for \ids to perform well.
Otherwise, only \dids can leverage multi-fidelity information to an advantage.

The current literature in the field of multi-fidelity blackbox optimization mainly studies the unconstrained case or model-based optimization approaches.
This research proposes methods to handle multi-fidelity constraints in the context of direct search methods.
Later work will incorporate the aforementioned objective function focused literature with the presented methods to consider the generalized case where all available multi-fidelity information is leveraged.

\section*{Acknowledgments}
The \solar{} problems are available at \href{https://github.com/bbopt/solar}{\url{https://github.com/bbopt/solar}}.
The NOMAD software package is available at \href{https://github.com/bbopt/nomad}{\url{https://github.com/bbopt/nomad}}.
This work is partly supported by the NSERC Alliance-Mitacs Accelerate grant ALLRP 571311-21 (``Optimization of future energy systems'') in collaboration with Hydro-Qu\'ebec. 
It is also funded by X.~Lebeuf's doctoral Hydro-Qu\'ebec excellence grant.

\bibliographystyle{plain}
\bibliography{bibliography}

\pdfbookmark[1]{References}{sec-refs}

\end{document}